\documentclass[11pt]{article}

\usepackage{amssymb}

\usepackage{amsmath}

\usepackage{amsfonts}

\usepackage{epsfig}

\usepackage{graphics}

\newtheorem{lemma}{Lemma}

\newtheorem{conjecture}[lemma]{Conjecture}

\begin{document}

\begin{titlepage}

\title{\bf{Random Matrix Theory and the Fourier Coefficients of
Half-Integral Weight Forms}}

\vspace {2 in}

\author{J.\ B.\ Conrey$^{1,2}$, J.\ P.\ Keating$^2$, M.\ O.\ Rubinstein$^3$ and
N.\ C.\ Snaith$^{2}$\\ \\
1.\ American Institute of Mathematics\\ 360 Portage Avenue\\Palo Alto, CA 94306\\USA\\ \\
2.\ School of Mathematics \\ University of Bristol\\ Bristol BS8 1TW\\ UK \\ \\
3.\ Pure Mathematics \\ University of Waterloo\\ 200 University Ave W\\Waterloo, ON, Canada\\N2L 3G1\\}

\date{\today}

\maketitle

\thispagestyle{empty}

\vspace{.5cm}

\begin{abstract}

Conjectured links between the distribution of values taken by the
characteristic polynomials of random orthogonal matrices and that
for certain families of $L$-functions at the centre of the
critical strip are used to motivate a series of conjectures
concerning the value-distribution of the Fourier coefficients of
half-integral weight modular forms related to these $L$-functions.
Our conjectures may be viewed
as being analogous to the Sato-Tate conjecture for integral weight
modular forms. Numerical evidence is presented in support of them.

\end{abstract}

\end{titlepage}

\maketitle

\section{Introduction}

\label{section:intro}

The limiting value distribution of the Fourier coefficients of
integral weight modular forms is given by the celebrated
Sato-Tate conjecture.
Our purpose here is to identify the
implications of some recent conjectures concerning the
value-distribution of certain families of $L$-functions at the
centre of the critical strip for the distribution of the Fourier
coefficients of related half-integral weight modular forms.

These conjectures concern the relationship between properties of
$L$-functions and random matrix theory.  It was conjectured by
Montgomery \cite{kn:mont73} that correlations between the zeros of
the Riemann zeta function on the scale of the mean zero separation
coincide with those between the phases of the eigenvalues of
unitary matrices, chosen at random, uniformly with respect to Haar
measure on the unitary group, in the limit of large matrix size.
This is supported both by theoretical \cite{kn:mont73, kn:rudsar,
kn:bogkea95, kn:bogkea96} and extensive numerical
\cite{kn:odlyzko89, kn:rub98} evidence.  Katz and Sarnak
\cite{kn:katzsarnak99a} then generalized the connection by
suggesting that statistical properties of the zeros within various
families of $L$-functions, computed by averaging over a given
family, coincide with those of the eigenvalues of random matrices
from the various classical compact groups,  the particular group
being determined by the family in question. Based on these ideas,
a link was proposed in \cite{kn:keasna00a} between the leading
order asymptotics of the value distribution of the Riemann zeta
function on its critical line and that of the characteristic
polynomials of random unitary matrices, giving, for example, 
an explicit conjecture for the leading order asymptotics of the
moments of the zeta function.  This approach was then extended to
relate the value distribution of $L$-functions, in families, at
the centre of the critical strip, to that of the characteristic
polynomials of matrices from the various classical compact groups
\cite{kn:confar00, kn:keasna00b}.  It has also recently been
extended to include all lower order terms in the asymptotics
\cite{kn:cfkrs}. (For other related results, see
\cite{kn:hughes00, kn:hughes01}.)  These developments have
recently been reviewed in \cite{kn:conrey1,kn:keasna03}.

Our strategy here is to combine these
random-matrix-theory-inspired conjectures for the value
distribution of $L$-functions with formulae due to Waldspurger
\cite{kn:waldspurger81}, Shimura \cite{kn:shimura73}, Kohnen
and Zagier \cite{kn:kohzag81}, and Baruch and Mao \cite{kn:barmao03}
which relate the values taken by
$L$-functions associated with elliptic curves at the centre of the
critical strip to the Fourier coefficients of certain
half-integral weight modular forms.  This approach was first
outlined in \cite{kn:ckrs00}, where its implications for the
vanishing of $L$-functions were examined.  Here we use the same
ideas in order to develop conjectures for the value-distribution
of the Fourier coefficients.

The outline of this paper is as follows.  In Section
\ref{section:modular} we give a brief overview of the relationship
between modular forms and $L$-functions associated with elliptic
curves. In Section \ref{section:rmt} we review the results from
Random Matrix Theory that we shall need, and some conjectures they
suggest for the value-distribution of $L$-functions associated
with elliptic curves. In Section \ref{section:coeffs}, we combine
results from Section \ref{section:modular} and Section
\ref{section:rmt} to motivate conjectures for the moments of the
Fourier coefficients, the value distribution of the logarithm of
the Fourier coefficients, and the distribution of the coefficients
themselves.  In Section \ref{section:vanishing} we review the
implications, first outlined in \cite{kn:ckrs00}, for the
vanishing of $L$-functions at the centre of the critical strip.
Finally, in Section \ref{section:numerics} we describe numerical
experiments whose results support our conjectures.

\section{Modular forms and $L$-functions}

\label{section:modular}

We review in this section the connection between the quadratic
twists of modular $L$-functions and the Fourier coefficients of
half-integral weight modular forms.  This is central to the
motivation underlying the conjectures we make in subsequent
sections.  In order to be concrete, we concentrate on the case of
quadratic twists of $L$-functions associated with elliptic
curves.

Let $L_E(s)$ be the $L$-function associated with an elliptic curve $E$ over
$\mathbb{Q}$ with Dirichlet series and Euler product given by
\begin{align}
\label{eq:9}
    L_{E}(s) & =\sum_{n=1}^\infty\frac{a_n}{n^s} 
    =\prod_{p\mid \Delta}
    \left(1-a_p p^{-s}\right)^{-1}
    \prod_{p\nmid \Delta}
    \left(1-a_p p^{-s}+p^{1-2s}\right)^{-1}\\
    & =
    \prod_{p} \mathcal{L}_p(1/p^s)
    , \quad \quad \Re(s) > 3/2,
\end{align}
with $\Delta$ the discriminant of $E$, and $a_p=p+1-\#E(\mathbb{F}_p)$, where
$\#E(\mathbb{F}_p)$ denotes the number of points on $E$ when regarded over ${F}_p$.
It is a consequence of the recently solved
Taniyama-Shimura conjecture~\cite{kn:wiles95}~\cite{kn:taywil95}~\cite{kn:BCDT01} 
that $L_E(s)$ has
analytic continuation to $\mathbb{C}$ and satisfies a
functional equation of the form
\begin{equation}
\label{eq:10a}
    \left(\frac{2\pi}{\sqrt{Q}}\right)^{-s}{\rm
    \Gamma}(s)L_E(s)=w_E\left(\frac{2\pi}{\sqrt{Q}}\right)^{s-2}{\rm\Gamma}(2-s)L_E(2-s),
\end{equation}
where $Q$ is the conductor of the elliptic curve $E$ and $w_E=\pm 1$.

We let $\chi_d(n)=\left(\frac{d}{n}\right)$
for a fundamental discriminant $d$, where
$\left(\frac{d}{n}\right)$ is the Kronecker symbol.
The twisted $L$-function
\begin{equation}
\label{eq:12}
    L_{E}(s,\chi_d)=\sum_{n=1}^{\infty}\frac{a_n\chi_d(n)}{n^s}
\end{equation}
is the $L$-function of the elliptic curve $E_d$, the quadratic
twist of $E$ by $d$.  
If ($d,Q)=1$, then $L_{E}(s,\chi_d)$ satisfies the functional equation
\begin{equation}
    \label{eq:14} \left(\frac{2\pi}{\sqrt{Q}|d|}\right)^{-s}
    {\rm \Gamma}(s)L_{E}(s,\chi_d)=\chi_d(-Q)w_E\left(\frac{2\pi}{\sqrt{Q}|d|}\right)^{s-2}
    {\rm \Gamma}(2-s)L_{E}(2-s,\chi_d).
\end{equation}
We shall be interested in the case when $w_E \chi_d(-Q)=1$, since 
otherwise $L_{E}(1,\chi_d)$ is trivially equal to zero.

We now come to a key point.  The $L$-functions above are
related to half-integer weight modular forms via
formulae due to Waldspurger \cite{kn:waldspurger81}, Shimura
\cite{kn:shimura73}, Kohnen and Zagier \cite{kn:kohzag81},
and Baruch and Mao \cite{kn:barmao03}. One must distinguish between
positive and negative $d$, and one must also sort them
according to various residue classes. This has been worked out
explicitly for thousands of examples by Rodriguez-Villegas and Tornaria
in the case that $Q$ is squarefree and they kindly supplied a database of 
such forms to the authors \cite{kn:rt04}.

Specifically, for $Q$ squarefree, assume that $d<0$ and that 
$\chi_d(p)=-a_p$ for every $p \mid Q$ ($a_p=\pm 1$ for $p \mid Q$ when $Q$ is squarefree).
Notice that such $d$ are restricted to
$\prod_{p\mid Q \atop p \text{odd}} ((p-1)/2)$ residue classes mod $Q$ or $4Q$,
depending on whether $Q$ is odd or even.
Then, for such $d$,
\begin{equation}
\label{eq:19} 
    L_E(1,\chi_d)=\kappa_E \frac{c_E(|d|)^2}{\sqrt{|d|}},
\end{equation}
where $c_E(|d|)\in \mathbb{Z}$ are the Fourier coefficients of a
weight $3/2$ modular form, and where $\kappa_E$ depends on $E$.

For $d>0$, the work of Baruch, Mao, Rodriguez-Villegas, and Tornaria
\cite{kn:barmao03} \cite{kn:mao-etal} gives the
relevant weight $3/2$ form for $Q$ prime and $d$'s satisfying $\chi_d(Q)=a_Q$.
One has the same relation as above, but with a different proportionality constant.
Some examples are listed in Section~\ref{section:numerics}.
In either case, given a coefficient of the weight $3/2$ form,
the constant $\kappa_E$ can be evaluated either by comparison with the
Birch and Swinnerton-Dyer conjecture, or by numerically computing 
$L_E(1,\chi_d)$ as a series involving the exponential function.

Our strategy will be to write down conjectures for the value
distribution of the Fourier coefficients $c(|d|)$ by coupling the
connection (\ref{eq:19}) with conjectures motivated by random
matrix theory for the value distribution of $L_{E}(1,\chi_d)$.

\section{$L$-functions and random matrices}

\label{section:rmt}

It was conjectured by Montgomery \cite{kn:mont73} that the zeros
of the Riemann zeta function are distributed statistically like
the eigenvalues of random hermitian (self-adjoint) matrices, or,
equivalently, like the phases of the eigenvalues of random unitary
matrices.  This extends to the zeros of any given principal
$L$-function \cite{kn:rudsar}.  It was conjectured by Katz and
Sarnak \cite{kn:katzsarnak99a} that the distribution of zeros
defined by averaging over families of $L$-functions with the
height up the critical line fixed coincides with the distribution
of the phases of the eigenvalues of matrices from one of the
classical compact groups, depending on the family in question.

These ideas motivated the conjecture \cite{kn:keasna00a} that, asymptotically,
the moments of the Riemann zeta function (or any other
principal $L$-function) averaged high on its critical line coincide, up to a simple
arithmetical factor, with the moments of the characteristic polynomials of random
unitary matrices.
This suggestion was then extended to relate the moments
of families of $L$-functions at the centre of the critical strip to those
of the characteristic polynomials of matrices from the various classical
compact groups \cite{kn:confar00,kn:keasna00b}.

For any elliptic curve $E$ it is conjectured that the family
of even-functional equation quadratic twists 
\begin{equation}
\label{eq:25}
 \Phi_E=\{L_E(s,\chi_d):w_E\chi_d(-Q)=+1\}
\end{equation}
is orthogonal.  Specifically, this family conjecturally has
symmetry type $O^{+}$.  Thus the value distribution of
$L_E(1,\chi_d)$ should be related to that of the characteristic
polynomials of matrices in $SO(2N)$, at the spectral symmetry
point, with $N \sim \log(|d|)$.

For an orthogonal matrix $A$, the characteristic polynomial may
be defined by
\begin{equation}
\label{eq:26}
Z(A,\theta)=\det\left( I-A{\rm e}^{-i\theta}\right).
\end{equation}
The eigenvalues of $A$ form complex conjugate pairs ${\rm e}^{\pm
i\theta_n}$, and so the symmetry point is at $\theta=0$.  The
moments of $Z(A,0)$ are defined by averaging over $A$ with respect
to normalized Haar measure for $SO(2N)$, $dA$:
\begin{equation}
\label{eq:27}
M_O(N,s)=\int_{SO(2N)}|\det\left( I-A\right) |^sdA.
\end{equation}
It was shown in \cite{kn:keasna00b} that for ${\rm Re}s>-1/2$
\begin{equation}
    \label{eq:28}
    M_O(N,s)=2^{2Ns}\prod_{j=1}^N\frac{{\rm
    \Gamma}(N+j-1){\rm \Gamma}(s+j-1/2)}{{\rm \Gamma}(j-1/2){\rm
    \Gamma}(s+j+N-1)}
\end{equation}
and that as $N\rightarrow \infty$
\begin{equation}
\label{eq:29}
 M_O(N,s)\sim g_s(O^+)N^{s(s-1)/2},
\end{equation}
with
\begin{equation}
\label{eq:30}
g_s(O^+)=\frac{2^{s^2/2}G(1+s)\sqrt{{\rm \Gamma}(1+2s)}}
{\sqrt{G(1+2s){\rm \Gamma}(1+s)}},
\end{equation}
where $G$ is Barnes' $G$-function:
\begin{equation}
\label{eq:31}
G(z+1)=(2\pi)^{z/2}\exp(-((\gamma+1)z^2+z)/2)\prod_{n=1}^{\infty}(1+\frac{z}{n})^n\exp(-z+z^2/2n).
\end{equation}
It follows from the fact that $G(1)=1$ and $G(z+1)={\rm \Gamma}(z)G(z)$ that
\begin{equation}
\label{eq:32}
g_k(O^+)=2^{k(k+1)/2}\prod_{j=1}^{k-1}\frac{j!}{2j!},
\end{equation}
for integer $k$.
Note that the right-hand side of (\ref{eq:28}) has a meromorphic continuation to the whole complex
$s$-plane.

For integer $k$, $M_O(N,k)$ is a polynomial of order $k(k-1)/2$:
\begin{equation}
M_O(N,k)=\left(2^{k(k+1)/2}
\prod_{j=1}^{k-1}\frac{j!}{(2j)!}\right) \prod_{0\leq i<j\leq k-1}
(N+\frac{i+j}{2}).
\end{equation}
It can also be written in terms of a multiple contour integral
\cite{kn:cfkrs1}, a form that will be useful later for comparison
with $L$-functions:
\begin{eqnarray}
\label{eq:contint}
 M_O(N,k)&=&\frac{(-1)^{k(k-1)/2}2^k}{(2\pi i)^k
k!} \oint \cdots
\oint e^{N\sum_{j=1}^kz_j} \\
&&\quad\times \prod_{1\leq \ell<m\leq k}
(1-e^{-z_m-z_{\ell}})^{-1} \frac{\Delta(z_1^2,\cdots,z_k^2)^2}
{\prod_{j=1}^k z_j^{2k-1}} dz_1\cdots dz_k.
\end{eqnarray}
Here $\Delta(x_1,\ldots,x_n)=\prod_{1\leq i<j\leq n}(x_j-x_i)$ and
the contours enclose the poles at zero.

The value distributions of the characteristic polynomial and its
logarithm can be written down directly using (\ref{eq:28}).  Let
$P_N(t)$ denote the probability density function associated with
the value distribution of $|Z(A,0)|$, i.e.
\begin{equation}
{\rm meas.}\{A\in SO(2N): \alpha<|{\rm det}(I-A)|\le
\beta\}=\int_{\alpha}^{\beta}P_N(t)dt. \label{eq:33}
\end{equation}
Then
\begin{equation}
P_N(t)=\frac{1}{2\pi
it}\int_{c-i\infty}^{c+i\infty}M_O(N,s)t^{-s}ds, \label{eq:34}
\end{equation}
for any $c>0$.

The asymptotics of $P_N(t)$ as $t\rightarrow 0$ comes from the
pole of $M_O(N,s)$ at $s=-1/2$.  The residue there is
\begin{equation}
h(N)=\frac{2^{-N}}{{\rm \Gamma}(N)}\prod_{j=1}^{N}\frac{{\rm \Gamma}(N+j-1){\rm \Gamma}(j)}
{{\rm \Gamma}(j-1/2){\rm \Gamma}(j+N-3/2)}
\label{eq:35}
\end{equation}
which, as $N \rightarrow \infty$, is given asymptotically by
\begin{equation}
h(N)\sim 2^{-7/8}G(1/2)\pi^{-1/4}N^{3/8}.
\label{eq:36}
\end{equation}
Thus
\begin{equation}
P_N(t)\sim h(N)t^{-1/2}
\label{eq:37}
\end{equation}
as $t\rightarrow 0$.

Importantly for us, one may deduce a central limit theorem for
$\log|Z(A,0)|$ from equation~(\ref{eq:34}) (see~\cite{kn:keasna00b}):
\begin{equation}
    \lim_{N\rightarrow \infty}{\rm meas.} \{A\in SO(2N):
    \frac{\log|{\rm det}(I-A)|+\frac{1}{2}\log N}{\sqrt{\log N}}\in
    (\alpha, \beta)\}
    =\frac{1}{\sqrt{2\pi}}\int_{\alpha}^{\beta}\exp(-\frac{t^2}{2})dt.
    \label{eq:38}
\end{equation}

The above results, proved for the characteristic polynomials of
random matrices, suggest the following conjectures for
$L_E(1,\chi_d)$ (these are special cases of those made in
\cite{kn:keasna00b}, but with the number theoretical
details worked out explicitly). 

Recall that we are assuming that $Q$ is squarefree.
Let $(d,Q)=1$. For $L_E(s,\chi_d)$ to have an even functional equation,
one needs $w_E \chi_d(-Q)=1$. This imposes a condition, in the case that
$Q$ is odd, on $d \mod Q$, and, in the case that $Q$ is even
on $d \mod 4Q$ ($4Q$ because $\chi_d(2)$ is periodic with period 8).
Let
\begin{equation}
\label{eq:Q tilde}
    \tilde{Q} = \begin{cases}
                    Q & \text{if $Q$ is odd and squarefree}\\
                    4Q & \text{if $Q$ is even and squarefree}.
                \end{cases}
\end{equation}
Next, we focus our attention on a subset of the $d$'s according to
certain residue classes mod $\tilde{Q}$.
We let
\begin{equation}
\label{eq:S-}
    S^-(X) = S^-_E(X) = \{ -X \leq d < 0; \chi_d(p) =  -a_p \ \ \text{for all $p \mid Q$}\}
\end{equation}
i.e. the set of negative fundamental discriminants $d$ up to $-X$, but restricted
according to a condition on $d$ mod $\tilde{Q}$.
Let
\begin{equation}
    \int_0^\alpha W^{-}_E(X,t) dt
    =\frac{\left| \{d \in S^-(X); L_E(1,\chi_d) < \alpha\} \right|}
                    {\left| S^-(X) \right|}
    \label{eq:40}
\end{equation}
and
\begin{equation}
    M^-_E(X,s)=\frac{1}
               {\left| S^-(X) \right|}
             \sum_{d \in S^-(X)} L_E(1,\chi_d)^s.
\label{eq:41}
\end{equation}
i.e. the $s$th moment of $L_E(1,\chi_d)$. 

For elliptic curves $E$ whose conductor $Q$ is prime, we
also look at the set of positive fundamental discriminants
\begin{equation}
\label{eq:S+}
    S^+(X) = S^+_E(X) = \{ 0 < d \leq X; \chi_d(Q) =  a_Q \}
\end{equation}
and define $W^+_E(X,t), M^+_E(X,s)$ as in the negative case.

Note that
\begin{equation}
    W^{\pm}_E(X,t)=\frac{1}{2\pi it}\int_{c-i\infty}^{c+i\infty}M^{\pm}_E(X,s)t^{-s}ds. \label{eq:42}
\end{equation}
The central conjecture is that as $X\rightarrow \infty$
\begin{equation}
M^{\pm}_E(X,s)\sim A^{\pm}(s)M_O(\log X, s) \label{eq:43}
\end{equation}
where
\begin{eqnarray}
    A^{\pm}(s)= &&\prod_{p \nmid Q} \left( 1-p^{-1}\right)^{s(s-1)/2} 
          \left(\frac{p}{p+1}\right)
          \left(
              \frac1p + \frac12 
              \left( 
                  \mathcal{L}_p(1/p)^{s} +
                  \mathcal{L}_p(-1/p)^{s} 
              \right)
          \right) \\ \notag
          \times&&\prod_{p \mid Q} \left( 1-p^{-1}\right)^{s(s-1)/2} 
          \mathcal{L}_p(\pm a_p/p)^{s} 
\label{eq:44}
\end{eqnarray}
is an arithmetical factor that depends on $E$.
A heuristic for $A^{\pm}(s)$ is given
in~\cite{kn:cfkrs}. The following conjectures are then motivated by this.

First we consider the moments of elliptic curve $L$-functions.

\begin{conjecture}
For $M^{\pm}_E(X,s)$ defined as in (\ref{eq:41}),
\begin{equation}
    \lim_{X\rightarrow \infty}\frac{M^{\pm}_E(X,s)}{(\log X)^{s(s-1)/2}}=
    A^{\pm}(s)g_s(O^+). \label{eq:45}
\end{equation}
\end{conjecture}

Second, for large $X$,  as $t\rightarrow 0$ we have
\begin{equation}
    W^{\pm}_E(X,t)\sim A^{\pm}(-\tfrac{1}{2})h(\log X) t^{-1/2}, \label{eq:46}
\end{equation}
and thus by scaling $t$ so that $t=O((\log X)^{-\gamma})$, for
$\gamma>1$, we are led to

\begin{conjecture} If $\;\log X f(\log X)\rightarrow 0$ as $X\rightarrow \infty$,
then
\begin{equation}
\label{eq:doublelimit} \lim_{X\rightarrow \infty}
\frac{\sqrt{f(\log X)}}{(\log X)^{\tfrac{3}{8}}}W^{\pm}_E(X,f(\log
X)y)=B y^{-1/2}
\end{equation}
 where
\begin{equation}
B=2^{-7/8}G(1/2)\pi^{-1/4}A^{\pm}(-1/2) \label{eq:48}
\end{equation}
\end{conjecture}

In addition, we have

\begin{conjecture}
\begin{eqnarray}
\label{eq:49} 
    &&\lim_{X\rightarrow \infty}
    \frac{1}{\left| S^{\pm}(X) \right|}
    \left| \{ d \in S^{\pm}(X); 
              \frac{\log L_E(1,\chi_d)+\frac{1}{2}\log\log|d|}{\sqrt{\log\log|d|}} 
              \in (\alpha, \beta)
           \} 
    \right| \\
    &&\qquad=\frac{1}{\sqrt{2\pi}}\int_{\alpha}^{\beta}\exp(-\frac{t^2}{2})dt.\nonumber
\end{eqnarray}
We take $\log L_E(1,\chi_d)+\frac{1}{2}\log\log|d|$ to lie outside the interval if $L_E(1,\chi_d)=0$
\end{conjecture}

The above conjecture is similar to central limit theorems for the
Riemann zeta function and other $L$-functions usually attributed
to Selberg \cite{kn:selberg46,kn:selberg46a}. An analogous conjecture is
made in \cite{kn:keasna00b} for quadratic Dirichlet $L$-functions.

Further, we have a conjecture, closely related to the preceding
one, for the distribution of the full range of values of
$L_E(1,\chi_d)$.

\begin{conjecture}
\begin{eqnarray}
\label{eq:distconj} 
    &&\lim_{X\rightarrow \infty}
    \frac{1}{\left| S^{\pm}(X) \right|}
    \left| \{ d \in S^{\pm}(X); 
              \alpha\leq (\sqrt{\log|d|}
              L_E(1,\chi_d))^{\tfrac{1}{\sqrt{\log \log|d|}}}\leq \beta
           \} 
    \right| \\
    &&\qquad =\frac{1}{\sqrt{2\pi}} \int_{\alpha}^{\beta}
    \frac{1}{t}\; e^{-\tfrac{1}{2} (\log t)^2} dt\nonumber
\end{eqnarray}
for fixed $0\leq\alpha\leq \beta$.
\end{conjecture}

Finally, it is discussed in \cite{kn:cfkrs} how (\ref{eq:contint})
leads to conjectures for mean values of $L$-functions, not just at
leading order as in (\ref{eq:45}), but including all terms in the
expansion down to the constant term.

\begin{conjecture}
\label{conj:big}
With $k$ an integer
\begin{equation}
\label{eq:momint} 
    M^{\pm}_E(X,k)
    =
    \frac{1}{X}
    \int_0^X 
    \Upsilon^{\pm}_k
    \left(\log(t)
    \right) dt
    +O(X^{-\tfrac{1}{2}+\epsilon})
\end{equation}
as $X\rightarrow\infty$, where $\Upsilon_k$ is a polynomial of degree
$k(k-1)/2$ given by the $k$-fold residue
\begin{equation}
    \Upsilon^{\pm}_k(x)=\frac{(-1)^{k(k-1)/2}2^{k}}{k!} \frac{1}{(2\pi i)^k}
    \oint \cdots \oint
    \frac{F^{\pm}_k(z_1,\ldots,z_k)\Delta(z_1^2,\ldots,z_k^2)^2}
    {\prod_{j=1}^k z_j^{2k-1}} e^{x\sum_{j=1}^kz_j}dz_1\ldots dz_k,
\end{equation}
where the contours above enclose the poles at zero and
\begin{equation}
    F^{\pm}_k(z_1,\ldots,z_k)=A^{\pm}_k(z_1,\ldots,z_k) \prod_{j=1}^k 
    \left(
        \frac{\Gamma(1+z_j)}{\Gamma(1-z_j)}
        \left(\frac{Q}{4\pi^2} \right)^{z_j}
    \right)^{\tfrac{1}{2}}
    \prod_{1\leq i<j\leq k} \zeta(1+z_i+z_j)
\end{equation}
and $A^{\pm}_k$, which depends on $E$, is the Euler product which is absolutely convergent for
$\sum_{j=1}^k |z_j|<1/2$,
\begin{equation}
    A^{\pm}_k(z_1,\dots,z_k) =
    \prod_p F_{k,p}(z_1,\ldots,z_k)
    \prod_{1\le i < j \le k}
    \left(1-\frac{1}{p^{1+z_i+z_j}}\right)
\end{equation}
with, for $p \nmid Q$,
\begin{equation}
    F_{k,p}  = 
    \left(1+\frac 1 p\right)^{-1}\left(\frac 1 p +\frac{1}{2}
    \left(\prod_{j=1}^k
    \mathcal{L}_p \left(\frac{1}{p^{1 +z_j}} \right)+
    \prod_{j=1}^k\mathcal{L}_p
    \left(\frac{-1}{p^{1 +z_j}}\right) \right)\right) . 
\end{equation}
and, for $p\mid Q$,
\begin{equation}
    F_{k,p}=
    \prod_{j=1}^k
    \mathcal{L}_p \left(\frac{\pm a_p}{p^{1 +z_j}} \right).
\end{equation}

\end{conjecture}

\section{Conjectures relating to the value distribution of the Fourier coefficients of half-integral weight forms}

\label{section:coeffs}

Our goal now is to use the conjectures listed at the end of the previous section for the
value distribution of the $L$-functions associated with elliptic curves to motivate conjectures
for the value distribution of the Fourier coefficients of half-integral weight forms.
These follow straightforwardly from the connection (\ref{eq:19}).

In each case, we let $c(|d|)$ refer to the coefficients of the corresponding weight 3/2
modular form, as in~(\ref{eq:19}), and $\kappa_E^{\pm}$ refer to the corresponding 
proportionality constant.

Our first conjecture then follows from (\ref{eq:49}):

\begin{conjecture}{\rm (central limit conjecture)}
\begin{eqnarray}
\label{eq:}
    &&\lim_{X\rightarrow \infty}
    \frac{1}{\left| S^{\pm}(X) \right|}
    \left| \{ d \in S^{\pm}(X);
              \frac{2\log c(|d|) -\frac12 \log|d|+\frac12\log\log|d|}{\sqrt{\log\log|d|}}
              \in (\alpha, \beta)
           \}
    \right| \\
    &&\qquad=\frac{1}{\sqrt{2\pi}}\int_{\alpha}^{\beta}\exp(-\frac{t^2}{2})dt.\nonumber
\end{eqnarray}
We take $2\log c(|d|) -\frac12 \log(|d|)+\frac12\log\log|d|$ to lie outside the interval if $c(|d|)=0$
\end{conjecture}

This leads directly to a conjecture for the appropriately
normalized distribution of the coefficients themselves, which is
analogous to (\ref{eq:distconj}).

\begin{conjecture}
\begin{eqnarray}
    &&\lim_{X\rightarrow \infty}
    \frac{1}{\left| S^{\pm}(X) \right|}
    \left| \{ d \in S^{\pm}(X); 
              \alpha\leq 
                  \left(\frac{\kappa_E^{\pm}\sqrt{\log |d|}}
                  {\sqrt{|d|}}  c(|d|)^2\right)^
                  {\tfrac{1}{\sqrt{\log \log |d|}}}\leq \beta
           \} 
    \right| \\
    &&\qquad =\frac{1}{\sqrt{2\pi}} \int_{\alpha}^{\beta}
    \frac{1}{t}\; e^{-\tfrac{1}{2} (\log t)^2} dt\nonumber
\end{eqnarray}
for fixed $0\leq\alpha\leq \beta$.
\end{conjecture}

Our third conjecture follows from (\ref{eq:19}) and (\ref{eq:45}):

\begin{conjecture}{\rm (moment conjecture)}
\begin{equation}
    \lim_{X\rightarrow \infty}\frac{1}{(\log X)^{s(s-1)/2}}
    \frac{1}{|S^{\pm}(X)|}
    \sum_{d \in S^{\pm}(X)}
    \frac{c(|d|)^{2s}}{|d|^{s/2}}=(\kappa^{\pm}_E)^{-s}A^{\pm}(s)g_s(O^+).
\end{equation}
\end{conjecture}

Further, we have the conjecture following from (\ref{eq:19}) and
(\ref{eq:doublelimit}).

\begin{conjecture} If $\;\log X f(\log X)\rightarrow 0$ as $X\rightarrow \infty$,
then
\begin{equation}
\lim_{X\rightarrow \infty} 
    \frac{\sqrt{\kappa_E^{\pm} f(\log X)}}{(\log
    X)^{\tfrac{3}{8}}}
    \frac{1}{|S^{\pm}(X)|}
    \left| \{
        d \in S^{\pm}(X);
        \frac{c(|d|)^2}{\sqrt{|d|}}<f(\log X)y
    \} \right| 
    =B y^{-1/2}
\end{equation}
 where
\begin{equation}
B=2^{-7/8}G(1/2)\pi^{-1/4}A^{\pm}(-1/2)
\end{equation}
\end{conjecture}

Lastly, we have the analogue of (\ref{eq:momint})

\begin{conjecture}
With $k$ an integer, and summing over fundamental discriminants,
\begin{equation}
    \frac{(\kappa^{\pm}_E)^k}{|S^{\pm}(X)|}
    \sum_{d\in S^{\pm}(X)}
    \frac{c(|d|)^{2k}}{|d|^{k/2}}
    =
    \frac{1}{X}
    \int_0^X 
    \Upsilon^{\pm}_k
    \left(\log(t)
    \right) dt
    +O(X^{-\tfrac{1}{2}+\epsilon})
\end{equation}
as $X\rightarrow\infty$, where $\Upsilon^{\pm}_k$ is the polynomial of
degree $k(k-1)/2$ given in Conjecture \ref{conj:big}.
\end{conjecture}

The numerical evidence that supports these conjectures is amassed
in Section \ref{section:numerics}.

\section{Frequency of vanishing of $L$-functions}
\label{section:vanishing}

Examining the frequency of $L_E(1,\chi_d)=0$ as $d$ varies, as well
as the order of the zeros, has particular significance in the
context of the conjecture of Birch and Swinnerton-Dyer which says
that $L_E(s)$ has a zero at $s=0$ with order exactly equal to the
rank of the elliptic curve $E$.  Random matrix theory appears to
have a role in predicting the frequency of such zeros.  We have
argued in \cite{kn:ckrs00} that since (\ref{eq:19}) implies a
discretisation of the values of $L_E(1,\chi_d)$ and
\begin{equation}
    \int_0^\alpha W^{\pm}_E(X,t) dt
\end{equation}
is the probability that $L_E(1,\chi_d)$ has a value of $\alpha$ or
smaller, then the combination of (\ref{eq:42}) and (\ref{eq:43})
suggest the following.

\begin{conjecture}
There is a constant $c^{\pm}_E \geq 0$ such that
\begin{equation}
    \label{eq:number vanishings}
    \frac
    {\sum_{{d \in S^{\pm}(X) \atop |d| \text{prime}}\atop L_E(1,\chi_d)=0} 1}
    {\sum_{d \in S^{\pm}(X) \atop |d| \text{prime}} 1}
    \sim c^{\pm}_E X^{-1/4}(\log X)^{3/8}.
\end{equation}
\end{conjecture}
This conjecture first appeared in~\cite{kn:ckrs00}, but was stated with $c^{\pm}_E>0$.
However, it became clear in preparing numerics for this paper that $c^{\pm}_E$
can equal zero, and an arithmetic explanation
has been given for one of the examples in our data by Delaunay~\cite{kn:delaunay04}.

Here the fundamental discriminants have been restricted to prime
values to avoid extra two divisibility issues placed on the
coefficients $c(|d|)$. The constant $c^{\pm}_E$ remains
somewhat mysterious. The random matrix model suggests that it should be proportional
to $A^{\pm}(-1/2) \sqrt{\kappa^{\pm}_E}$. However, when one attempts to
apply the random matrix model to the problem of the discrete values taken
on by the $c(|d|)$ one ignores subtle arithmetic. It seems, from numerical experiments, 
that one needs to take into account further correction factors that
depend on the size of the torsion subgroup of $E$, but this is still not
understood. See Section 6.

Let $q \nmid Q$ be a fixed prime.
Another conjecture that follows from this approach concerns
sorting the $d$'s for which $L_E(1,\chi_d)=0$ according to residue classes
mod $q$, according to whether $\chi_d(q)=1$ or $-1$. Let
\begin{equation}
    R^{\pm}_q(X)=
    \frac{
        \sum_{{d \in S^{\pm}(X) \atop L_E(1,\chi_d)=0} \atop \chi_d(q)=1} 1
    }
    {
        \sum_{{d \in S^{\pm}(X) \atop L_E(1,\chi_d)=0} \atop \chi_d(q)=-1} 1
    }.
\end{equation}
\begin{conjecture}{\cite{kn:ckrs00}}
Let
\begin{equation}
    R_q = \left( \frac{q+1-a_q}{q+1+a_q} \right)^{1/2}.
\end{equation}
Then, for $q\nmid Q$,
\begin{equation}
    \lim_{X\rightarrow \infty}
    R^{\pm}_q(X)=R_q.
\end{equation}
\end{conjecture}
We believe this conjecture to hold even if we allow $d$ to range over 
different sets of discriminants, such as $|d|$ restricted to primes
(though in the latter case we must be sure to rule out there
being no vanishings at all due to arithmetic reasons).

A more precise conjecture given in~\cite{kn:crw04} incorporates the next term.
The lower terms do depend on whether we are looking at $S^{+}(X)$ as opposed
to $S^{-}(X)$.
We require some notation. Let $p$ be prime. For $p \mid Q$ set
\begin{equation}
    \beta(p) = \beta^{\pm}(p) =
    \begin{cases}
        \log(p)/(1+p) \quad \text{in the $+$ case} \\
        \log(p)/(1-p) \quad \text{in the $-$ case}
    \end{cases}
\end{equation}
and for $p \nmid Q$ set
\begin{equation}
    \beta(p) = \log(p)
    \left(
        \frac{(2-a_p)f_1(p)^{-1/2} + (2+a_p)f_2(p)^{-1/2}}
             {2+p\left(f_1(p)^{1/2} + f_2(p)^{1/2} \right)}
    \right)
\end{equation}
where
\begin{eqnarray}
    f_1(p)=1-a_p/p+1/p \notag \\
    f_2(p)=1+a_p/p+1/p \notag.
\end{eqnarray}
Next, let
\begin{equation}
    \lambda_{\nu}(q) = 
    \frac{\log(q)(\nu a_q -2)}{\nu a_q -q -1} -
    \frac{3\log(q)}{2(q-1)}
    -5 \gamma/2 + 
    \sum_{p\neq q} \left( \beta(p) -\frac{3\log(p)}{2(p-1)} \right)
\end{equation}
where $\gamma$ is Euler's constant.

\begin{conjecture}{\cite{kn:crw04}}
For $q\nmid Q$
\begin{equation}
    \label{eq:second term}
    R^{\pm}_q(X)=\left( \frac{q+1-a_q}{q+1+a_q} \right)^{1/2}
    \frac{g+\lambda_1(q)}{g+\lambda_{-1}(q)}
    +O(1/\log(X)^2)
\end{equation}
where
\begin{equation}
    g = \frac{8}{3}\log(XQ^{1/2}/(2\pi) )-1
\end{equation}
\end{conjecture}
These conjectures are supported by numerical evidence that will be
described in the next section.

\section{Numerical Experiments}
\label{section:numerics}

We present numerics for 2398 elliptic curves and millions of 
quadratic twists for each curve ($|d| < 10^8$)
confirming the aforementioned conjectures.
To test these conjectures we used the relation~(\ref{eq:19}).
To make our examples explicit, we list in Table~\ref{tab:ternary1} relevant data for
26 of the 2398 elliptic curves examined. The remaining data may be obtained
from the $L$-function database of one of the authors~\cite{kn:rub04}.
We used as the starting point for our computations a
database of Tornaria and Rodriguez-Villegas~\cite{kn:rt04}
that lists, for thousands of elliptic curves, the relevant 
ternary quadratic forms.

Each entry in this table contains the following data
\begin{table}[h]
\centerline{\small
\begin{tabular}{|l|l|l|l|l|}
\hline
name & $[a_1, a_2, a_3, a_4, a_6]$ & $\kappa$ & Number of $d$ & number ternary forms  \\ \hline
\multicolumn{5}{|l|}{relevant residue classes $|d|$ mod $\tilde{Q}$} \\ \hline
\multicolumn{5}{|l|}{linear combination} \\ \hline
\multicolumn{5}{|l|}{ternary forms} \\ \hline
\end{tabular}
}
\label{tab:ternarysample}
\end{table}

The name we use for an elliptic curve $E$ corresponds to Cremona's table
\cite{kn:cremona}, but with an extra subscript, either `i' or `r', standing for
imaginary or real respectively, to specify
whether we are looking at quadratic twists $L_E(s,\chi_d)$ with $d<0$ or $d>0$.
The naming convention used by Cremona includes the conductor $Q$ in the name for
the elliptic curve. So, for example, the first entry has name $11A_i$ which is the
elliptic curve of conductor 11. The `i' indicates that we are looking at
quadratic twists of $11A$ with $d<0$.

$[a_1, a_2, a_3, a_4, a_6]$ refers to the coefficients defining the equation of $E$,
\begin{equation}
    y^2 + a_1 xy + a_3 y = x^3 + a_2 x^2 + a_4 x + a_6.
\end{equation}
For each curve, our data consists of values of $L_E(1,\chi_d)$ with $|d|<10^8$, and $d$
restricted to $\prod_{p\mid Q \atop p \text{odd}} ((p-1)/2)$ residue classes mod $Q$ or $4Q$,
depending on whether $Q$ is odd or even. For imaginary twists, all the
Tornaria-Rodriguez examples have $Q$ squarefree, and in the case of real twists, all 
their examples have $Q$ prime.
The relevant residue classes and relevant modulus
is listed in the second line of each entry, and the number of such $d$ up to $10^8$ is given
by the `Number of $d$'.

Our $L(1,\chi_d)$ values are expressed in terms of the coefficients of a weight 3/2 modular form of 
level $4Q$ in the case of imaginary quadratic twists, $d<0$, and of
level $4Q^2$ in the case of real quadratic twists, $d>0$.
This weight 3/2 form is expressed as a linear combination of theta series attached
to positive definite ternary quadratic forms. The number of forms, say $r$,
is given in the last part of the first line,
while the linear combination, $[\alpha_1,\ldots,\alpha_r]$,
and ternary forms occupy the last two lines of each entry.
Each ternary form is specified as a sextuple of integers $\beta=[\beta_1,\beta_2,\beta_3,\beta_4,\beta_5,\beta_6]$.
The ternary form is
\begin{equation}
    f_\beta(x,y,z) = \beta_1 x^2 +\beta_2 y^2 + \beta_3 z^2 + \beta_4 yz + \beta_5 xz + \beta_6 xy
\end{equation}
and the theta associated series is given by
\begin{equation}
    \theta_\beta(w) =  \sum_{(x,y,z) \in \mathbb{Z}^3} w^{f_\beta(x,y,z)}
\end{equation}
Given this data, one defines
\begin{equation}
    \sum_{n=1}^\infty c(n) w^n = \sum_{j=1}^r \alpha_j \theta_{\beta_j}(w).
\end{equation}
Then, for fundamental discriminants lying in the relevant residue classes in the 
table (and $d<0$ or $d>1$ according to the name of the entry),
one has
\begin{equation}
    L_E(1,\chi_d) = \kappa c(|d|)^2 / |d|^{1/2}.
\end{equation}
Note, our values of $\kappa$ and $c(|d|)$ differ slightly from the values given
in \cite{kn:mao-etal} in that, in our tables, we ignored the value of
$c(1)$, for example in the case that $d>0$, and normalized the $\alpha_j$'s so that the 
gcd of all the $c(n)$'s, for $n \neq 1$, is equal to one, absorbing if necesary an extra 
square factor into $\kappa$. For example, while \cite{kn:mao-etal} gives for the curve $11A_r$
a value of $\kappa=.25384186\ldots$, we list a value of $\kappa$ that is $5^2=25$
times as big, but give $c(n)$'s which are $1/5$ as large.

The most comprehensive test carried out \cite{kn:cfkrs}
for moments involved looking
at conjecture 5 for millions of values of $L_{E_{11A}}(1,\chi_d)$,
with $d<0$ and $|d|=1,3,4,5,9 \mod 11$.
The coefficients, from \cite{kn:cfkrs}, of the polynomials $\Upsilon^{-}_k$ in 
conjecture 5 are given by Table~\ref{tab:fr}. In checking conjecture 5, we 
compared 
\begin{equation}
     \sum_{ {-850000<d<0}\atop{|d| =  1,3,4,5,9 \mod 11} }
      L_{11A_i}(\frac12 ,\chi_d)^k
\label{eqn:L11sum}
\end{equation}
to
\begin{equation}
     \sum_{ {-850000<d<0}\atop{|d| =  1,3,4,5,9 \mod 11} }
     \Upsilon^{-}_k\left(\log(|d|)\right).
\label{eqn:L11sum2}
\end{equation}
This comparison is depicted in Table \ref{tab:L11moments}.

\begin{table}[h!tb]
\centerline{\small
\begin{tabular}{|c|c|c|c|c|}
\hline
r& $f_r(1)$ &$f_r(2)$ &$f_r(3)$ &$f_r(4)$ \cr \hline
0 & 1.2353 & .3834 & .00804 & .0000058 \cr
1 & & 1.850 & .209 & .000444 \cr
2 & & & 1.57 & .0132 \cr
3 & & & 2.85 & .1919 \cr
4 & & & & 1.381 \cr
5 & & & & 4.41 \cr
6 & & & & 4.3 \cr
\hline
\end{tabular}
}
\caption{
Coefficients of $\Upsilon^{-}_k(x) = f_0(k) x^{k(k-1)/2} + f_1(k)
x^{k(k-1)/2-1} + \ldots$, for k=1,2,3,4.
}\label{tab:fr}\end{table}

\begin{table}[h!tb]
\centerline{\small
\begin{tabular}{|c|c|c|c|}
\hline
$k$ & (\ref{eqn:L11sum}) & (\ref{eqn:L11sum2})& ratio \cr \hline
1 &14628043.5       &14628305.       &0.99998 \cr 2 &100242348.8
&100263216.      &0.9998 \cr 3 &1584067116.8     &1587623419.
&0.998 \cr 4 &41674900434.9    &41989559937.    &0.993 \cr \hline
\end{tabular}
}
\caption{ Moments of $L_{11A_i}(1 ,\chi_d)$ versus
their conjectured values, for fundamental discriminants
$-85,000,000 < d < 0$, $|d|=1,3,4,5,9 \mod 11$, and $k=1,\ldots,4$.
The data agree with our conjecture to the accuracy to which we
have computed the moment polynomials $\Upsilon^{-}_k$.
}\label{tab:L11moments}\end{table}

The first part of Conjecture 1 is, for integral moments, a weaker form
of conjecture 5, while conjectures 8 and 10 follow from 1 and 5.

\begin{figure}[htp]
    \centerline{
         \psfig{figure=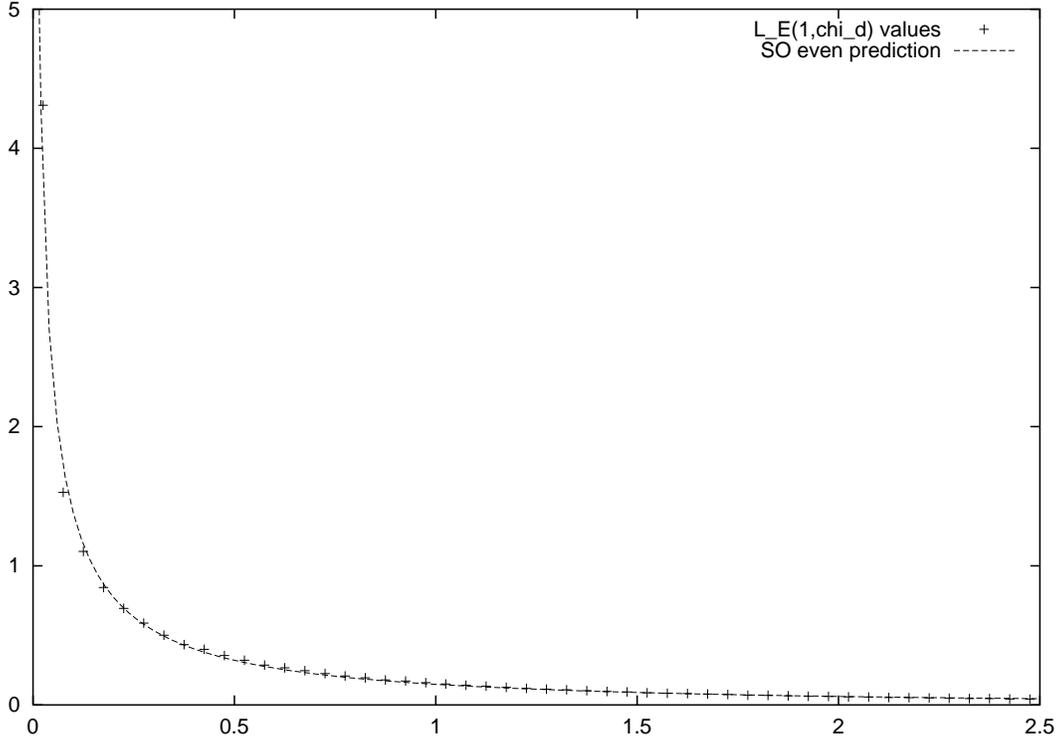,width=4in,angle=-90}
    }
    \caption{
     Value distribution 
     of $L_{E_{11}}(1,\chi_d)$, with $-85000000<d<0$, $d=2,6,7,8,10 \mod 11$.}
     \label{fig:value dist}
\end{figure}

In Figure~\ref{fig:value dist} we depict the numerical value distribution
of $L_{11A_i}(1 ,\chi_d)$, for fundamental discriminants 
$-85,000,000 < d < 0$, $|d|=1,3,4,5,9 \mod 11$,  compared to 
$P_N(t)$ obtained by taking the inverse Mellin transform of~(\ref{eq:28})
with $N=20$.
Because we are neglecting the arithmetic factor in computing the density,
a slight cheat was used to get a better fit. 
The histograms were rescaled by a constant along both axis until the histogram 
displayed matched
up nicely with the solid curve. Considering that we are compensating
for leaving out the arithemtic factor in such a naive way, it is a bit surprising how
nicely the two fit. The main point we wish to make is that the histogram does
exhibit $t^{-1/2}$ behaviour near the origin in support of part 2 of Conjecture 1.

\begin{figure}[htp]
    \centerline{
            \psfig{figure=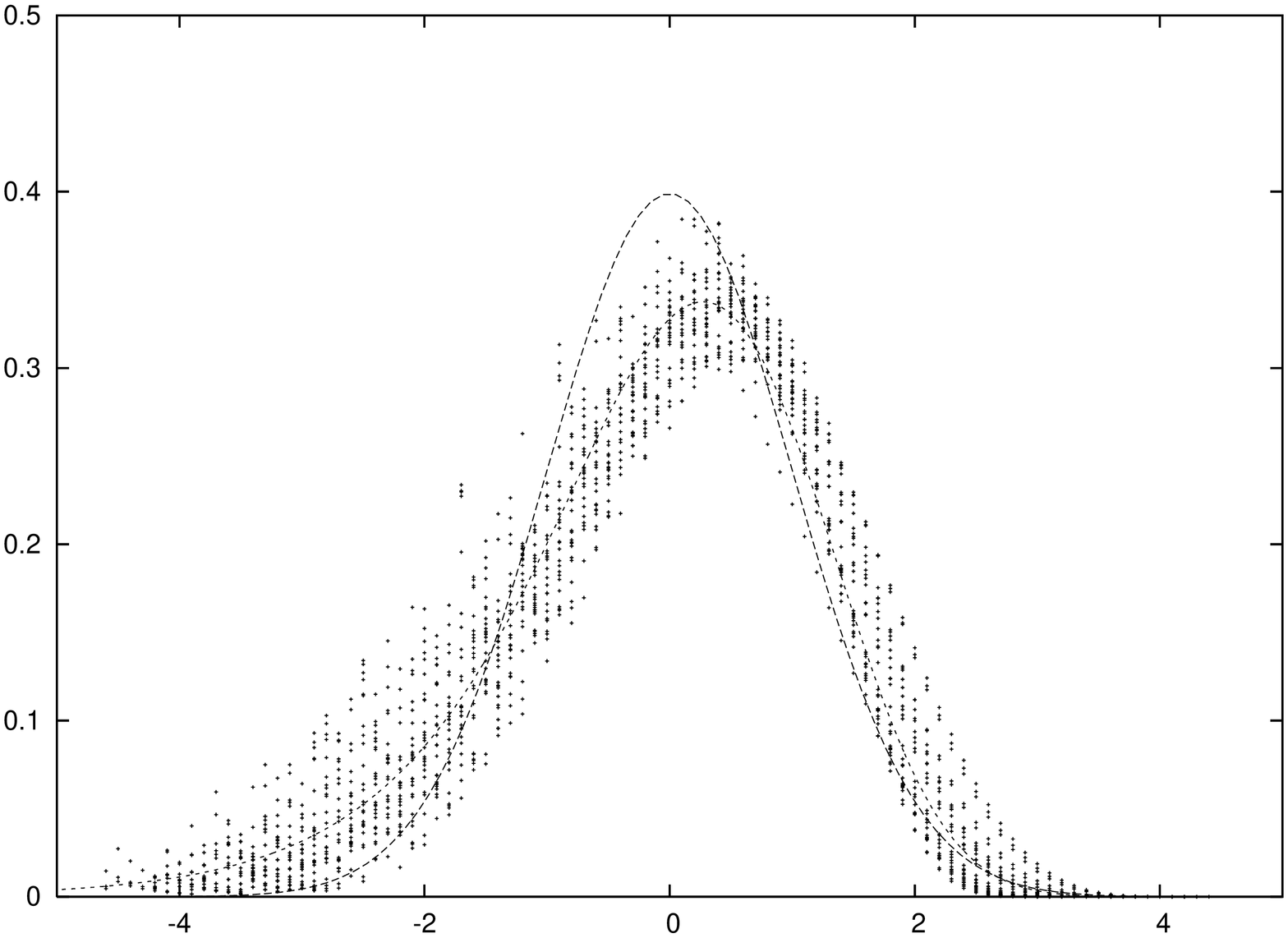,width=4in,angle=-90}
    }
    \centerline{
            \psfig{figure=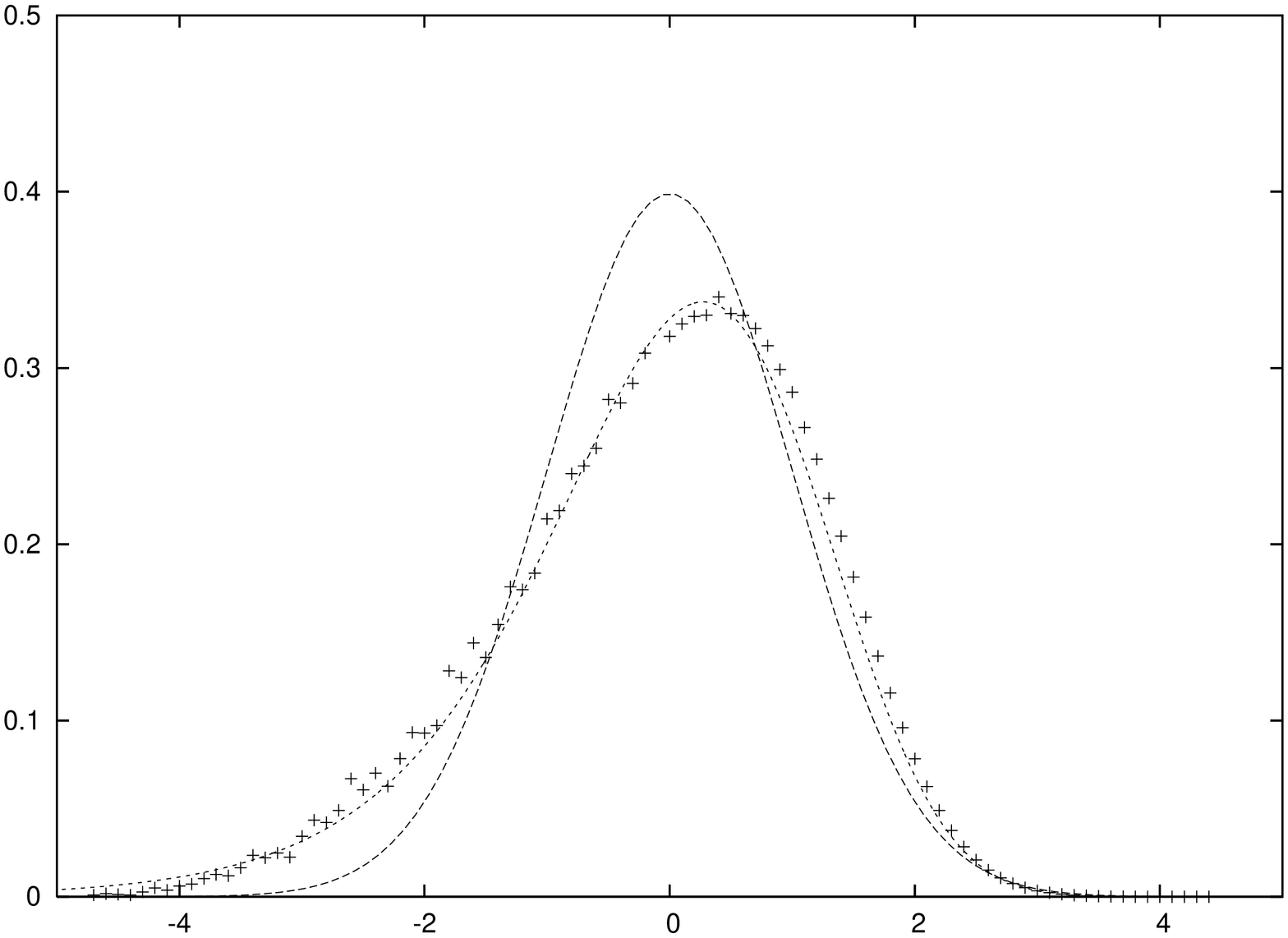,width=4in,angle=-90}
    }
    \caption{
        Value distribution of $\left( \log L_E(1,\chi_d)+\frac{1}{2}\log\log|d|\right)/\sqrt{\log\log|d|}$ 
        compared to the random matrix theory counterpart, $P_N(g(t)) g'(t)$ with $N=20$
        (rescaled as explained in the text), and
        its limit the standard Gaussian. In the first picture we superimpose the value distributions for
        the 26 curves described in Table~\ref{tab:ternary1}.
        The second picture shows the average value distribution of the 26.
    }
    \label{fig:gaussian}
\end{figure}

In Figure~\ref{fig:gaussian} we verify the central limit theorem described in
Conjecture~3. The first picture shows the value distributions, superimposed, of
\begin{equation}
    \frac{\log L_E(1,\chi_d)+\frac{1}{2}\log\log|d|}{\sqrt{\log\log|d|}}, \quad |d|<10^8
\end{equation}
for the 26 curves described in
Tables~\ref{tab:ternary1}. The second plot depicts the average value distribution
of the 26. These are compared against the standard Gaussian predicted in the
conjecture and also against the density function associated to the value distribution of
\begin{equation}
    \frac{\log |Z(A,0)| +\frac{1}{2}\log N}{\sqrt{\log N}}
\end{equation}
with $N=20$.
This density is given by
\begin{equation}
    P_N(g(t)) g'(t)
\end{equation}
where
\begin{equation}
    g(t)= \frac{\exp(t\sqrt{\log N}) }{N^{1/2}},
\end{equation}
and is shown in~\cite{kn:keasna00b} to tend to the standard Gaussian as $N\to \infty$.
To get a better fit to the numeric value distribution, one would also need to incorporate
the arithmetic factor. In the limit, this factor has no effect, but the convergence to the limit is
extremely slow. The variability of the arithmetic factor explains why the 26
densities superimposed in Figure~\ref{fig:gaussian} don't fall exactly on
top of one another. To get a better fit to the average, we set
$p(t)= P_N(g(t)) g'(t)$ and plot instead $\alpha p(\alpha t)$, with $\alpha=1.21$ chosen so 
that the density function visually lines up with the data.

\begin{figure}[htp]
    \centerline{
            \psfig{figure=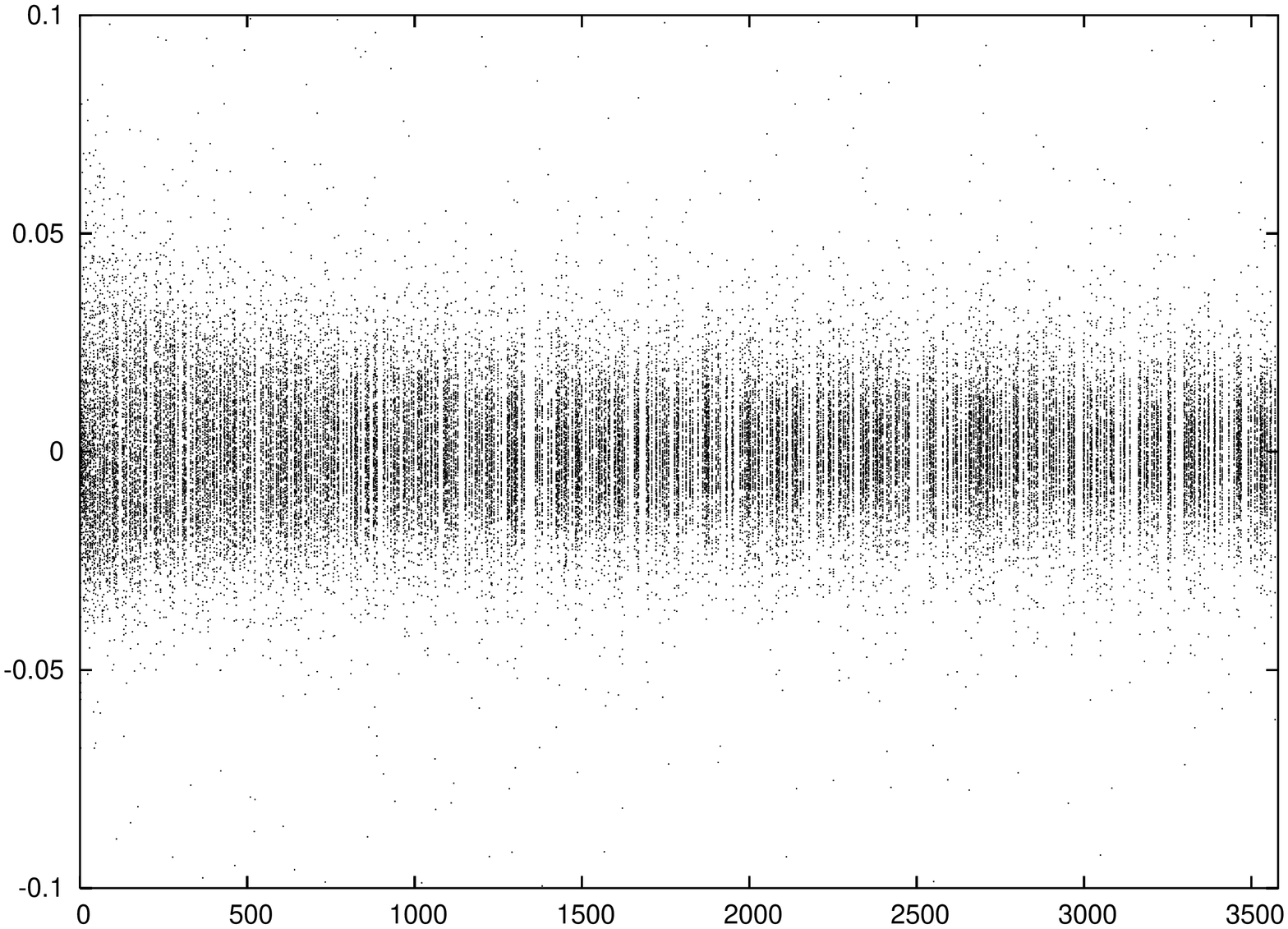,width=4in,angle=-90}
    }
    \centerline{
            \psfig{figure=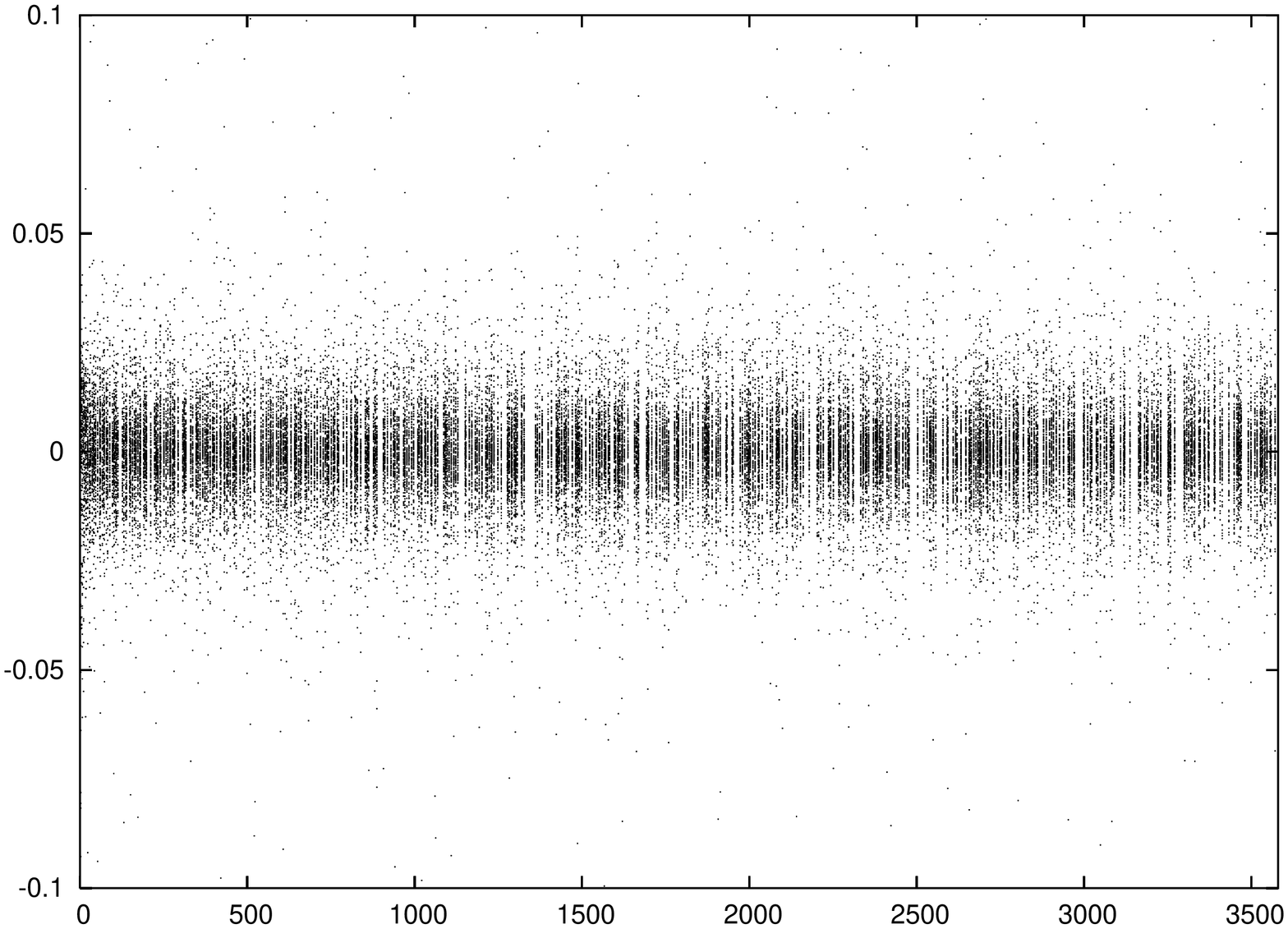,width=4in,angle=-90}
    }
    \caption
    {A plot \cite{kn:crw04} for one hundered elliptic curves of $R^{\pm}_q(X)-R_q$, top plot, and of (\ref{eq:conj 13})
    , bottom plot, for $2\leq q\leq 3571$, $X=10^8$.}
    \label{fig:1st vs 2nd}
\end{figure}

Conjecture 12 is verified numerically in the top plot in Figure~\ref{fig:1st vs 2nd}
which compares, for the first one hundred elliptic curves $E$ in our
database, and the sets $S^{\pm}_E(X)$,
the predicted value of $R_q$ to the actual value $R^{\pm}_q(X)$, with $X=10^8$.

The horizontal axis is $q$. For each $q$ and each of the one hundred
elliptic curves $E$ we plot $R^{\pm}_q(X)-R_q$, with $X=10^8$, and $q \leq 3571$.
For each $q$ on the horizontal there are 100 points corresponding to the 
100 values, one for each elliptic curve, of $R^{\pm}_q(X)-R_q$.
We see the values fluctuating
about zero, most of the time agreeing to within about $.02$.
The convergence in $X$ is predicted from secondary terms to be
logarithmically slow and one gets a better fit by including more terms
as predicted in conjecture 13.

This is depicted in the second plot of Figure~\ref{fig:1st vs 2nd}
which shows the difference
\begin{equation}
    \label{eq:conj 13}
    R^{\pm}_q(X)-
    \left( \frac{q+1-a_q}{q+1+a_q} \right)^{1/2}
    \frac{g+\lambda_1(q)}{g+\lambda_{-1}(q)}
\end{equation}
again with $X=10^8$, $q \leq 3571$, and the same one hundred
elliptic curves $E$. We see an improvement
to the first plot which uses just the main term.

\begin{figure}[htp]
    \centerline{
            \psfig{figure=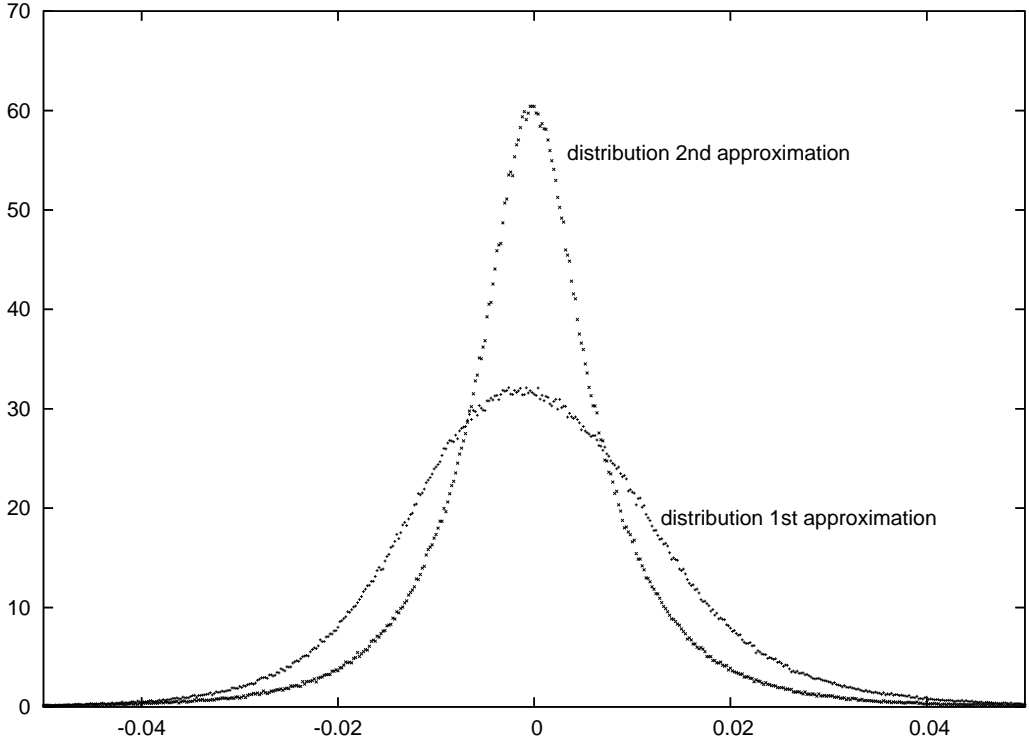,width=4in,angle=-90}
    }
    \caption
    {Distribution first approximation v.s. second approximation for ratio of vanishings \cite{kn:crw04}}
    \label{fig:dist 1st vs 2nd}
\end{figure}

This improvement is emphasized in Figure~\ref{fig:dist 1st vs 2nd}
which compares the distribution of $R^{\pm}_q(X)-R_q$ for all our $2398$ elliptic curves,
$X=10^8$, $q \leq 3571$, versus the distribution of~(\ref{eq:conj 13}).
The latter has smaller variance. These distributions
are not Gaussian. They depict the remainder of $R^{\pm}_q(X)$ compared to the first and second
conjectured approximations. There are yet further lower terms and these are given
by complicated sums involving the Dirichlet coefficients of $L_E(s)$.

\begin{figure}[htp]
    \centerline{
            \psfig{figure=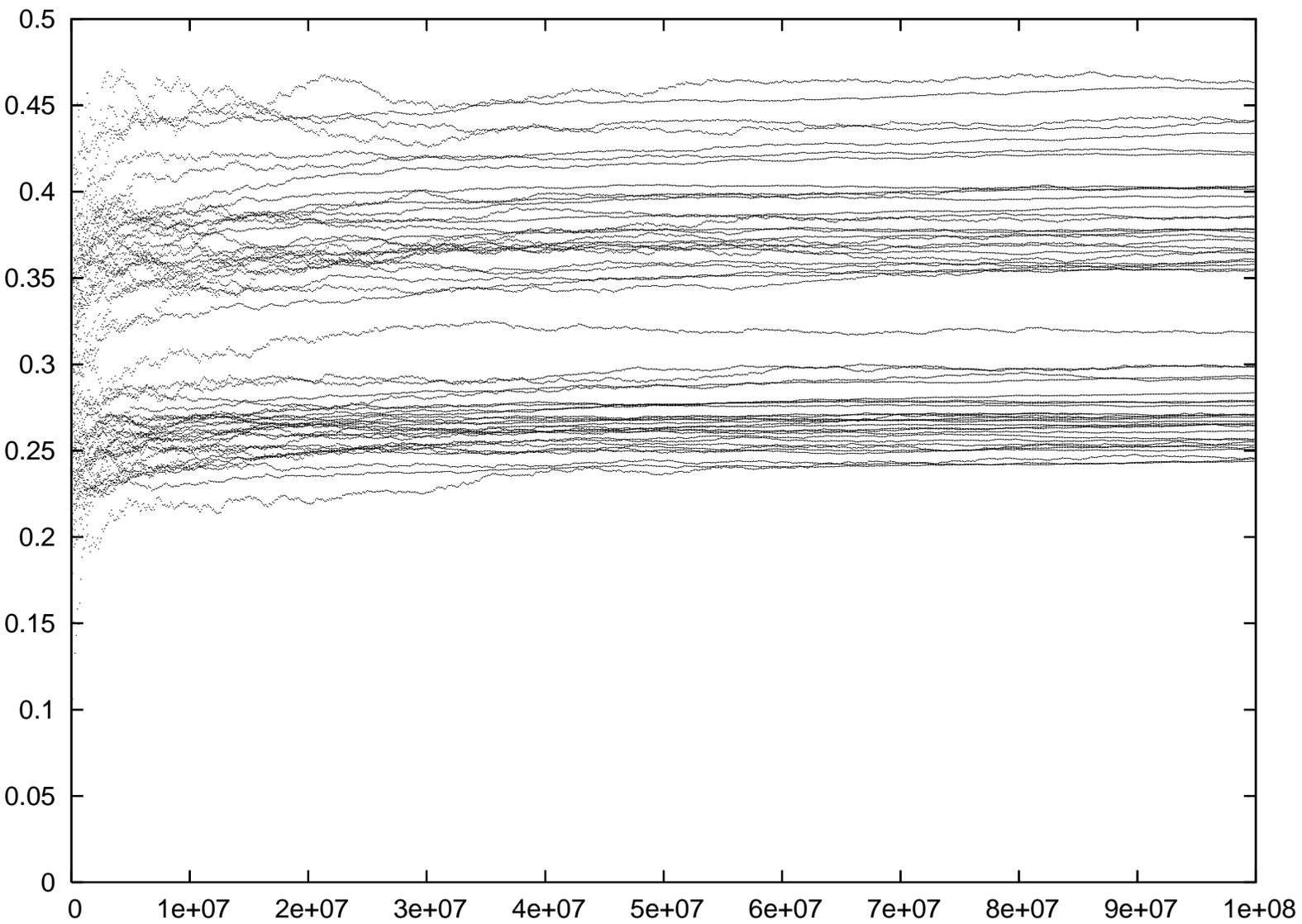,width=6in,angle=0}
    }
    \caption
    {A test of conjecture 11 for 55 elliptic curves in our database:
     we are plotting for each curve the l.h.s. of (\ref{eq:number vanishings})
     divided by $A^{\pm}(-1/2) \sqrt{\kappa^{\pm}_E} X^{-1/4}(\log X)^{3/8}$,
     $X=100000,200000,\ldots 10^8$. The graphs appear relatively flat.}
    \label{fig:c_E 76 curves}
\end{figure}

\begin{figure}[htp]
    \centerline{
            \psfig{figure=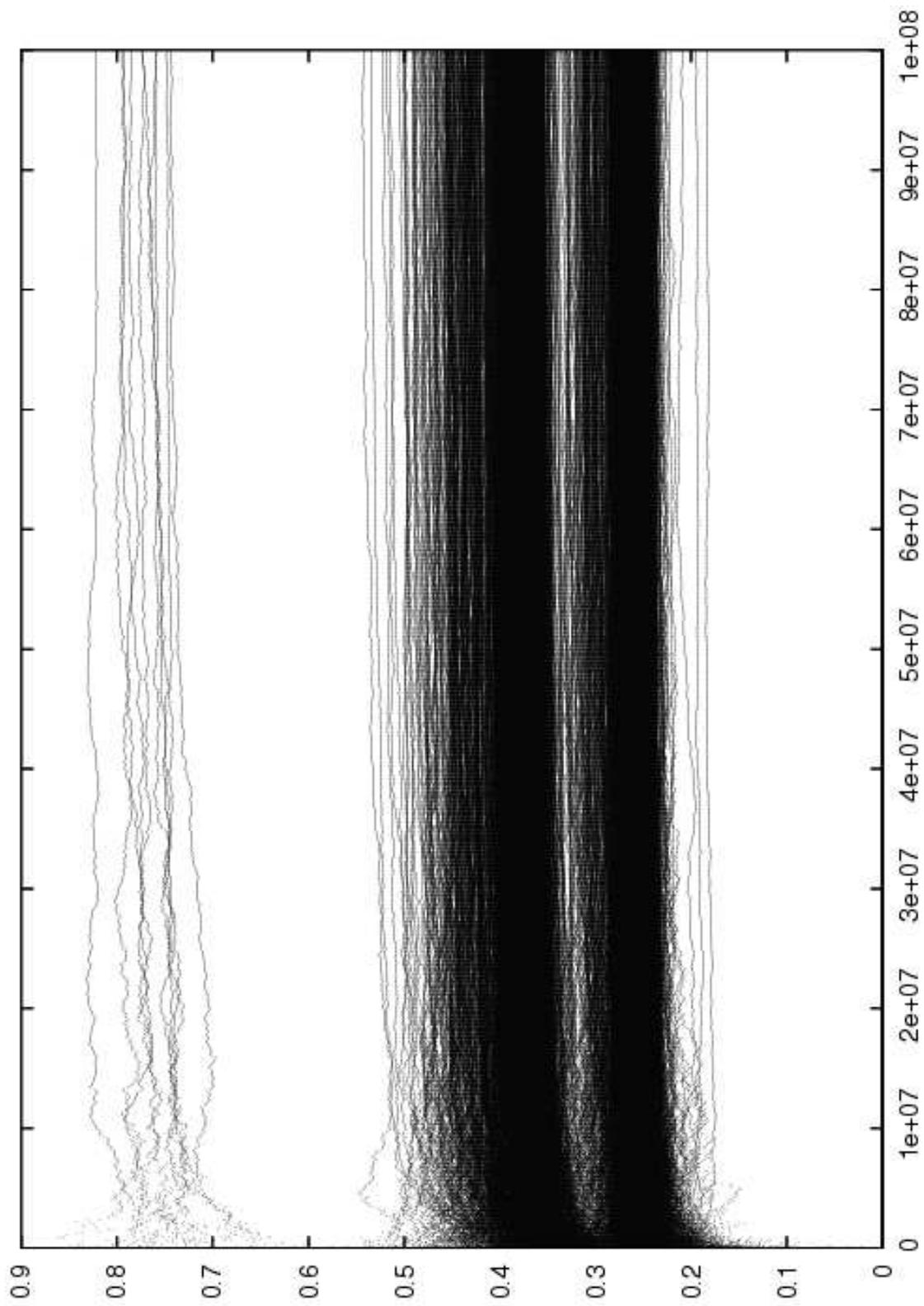,width=6in,angle=0} 
    }   
    \caption
    {The same as the previous figure, but for all elliptic curves in our database.}
    \label{fig:c_E all curves}
\end{figure}

To verify conjecture 11 we computed the l.h.s. of (\ref{eq:number vanishings})
divided by 
\begin{equation}
    A^{\pm}(-1/2) \sqrt{\kappa^{\pm}_E} X^{-1/4}(\log X)^{3/8}
\end{equation}
for our 2398 curves (545 of these have $c^{\pm}_E = 0$ and these are omitted) and
$X=100000,200000,\ldots 10^8$. The resulting values are depicted
in Figure \ref{fig:c_E all curves}. We also display in Figure \ref{fig:c_E 76 curves} 
the same data for a subset of 55 of these curves (those in our database with  $c^{\pm}_E \neq 0$
and whose conductor have leading digits $11$).
The graphs displayed in these two figures appear relatively flat.

\begin{figure}[htp]
    \centerline{
            \psfig{figure=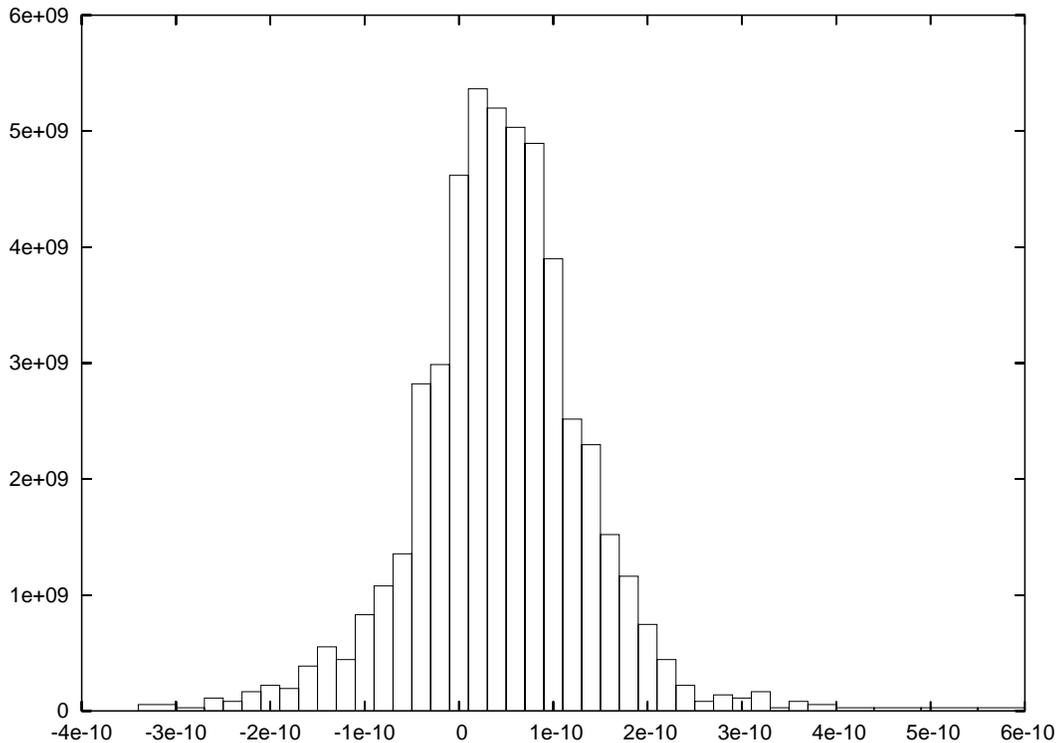,width=4in,angle=-90}
    }
    \caption
    {Distribution of the slopes of the graphs in Figure \ref{fig:c_E all curves}
     from $X=5\times 10^7$ to $X=10^8$. This tells us that the graphs are relatively flat.}
    \label{fig:slopes}
\end{figure}

To measure how flat these graphs are, we show in Figure \ref{fig:slopes}
the distribution of the slopes of the graphs in Figure \ref{fig:c_E all curves},
measured by sampling
each one at $X=5\times 10^7$ and $X=10^8$. Most have slopes that
are of size less than $10^{-10}$, and this suggests that the power of
$X$, namely $3/4$, in our conjecture is not off by more that $.01$, and
that the power of $\log(X)$, namely $3/8$ is not off by more than $.1$.
The mean does appear slightly to the right of 0, occuring at
$.5 \times 10^{-10}$, but this might be the result of using
a limited number of curves and also perhaps due to lower
order terms.

\begin{figure}[htp]
    \centerline{
            \psfig{figure=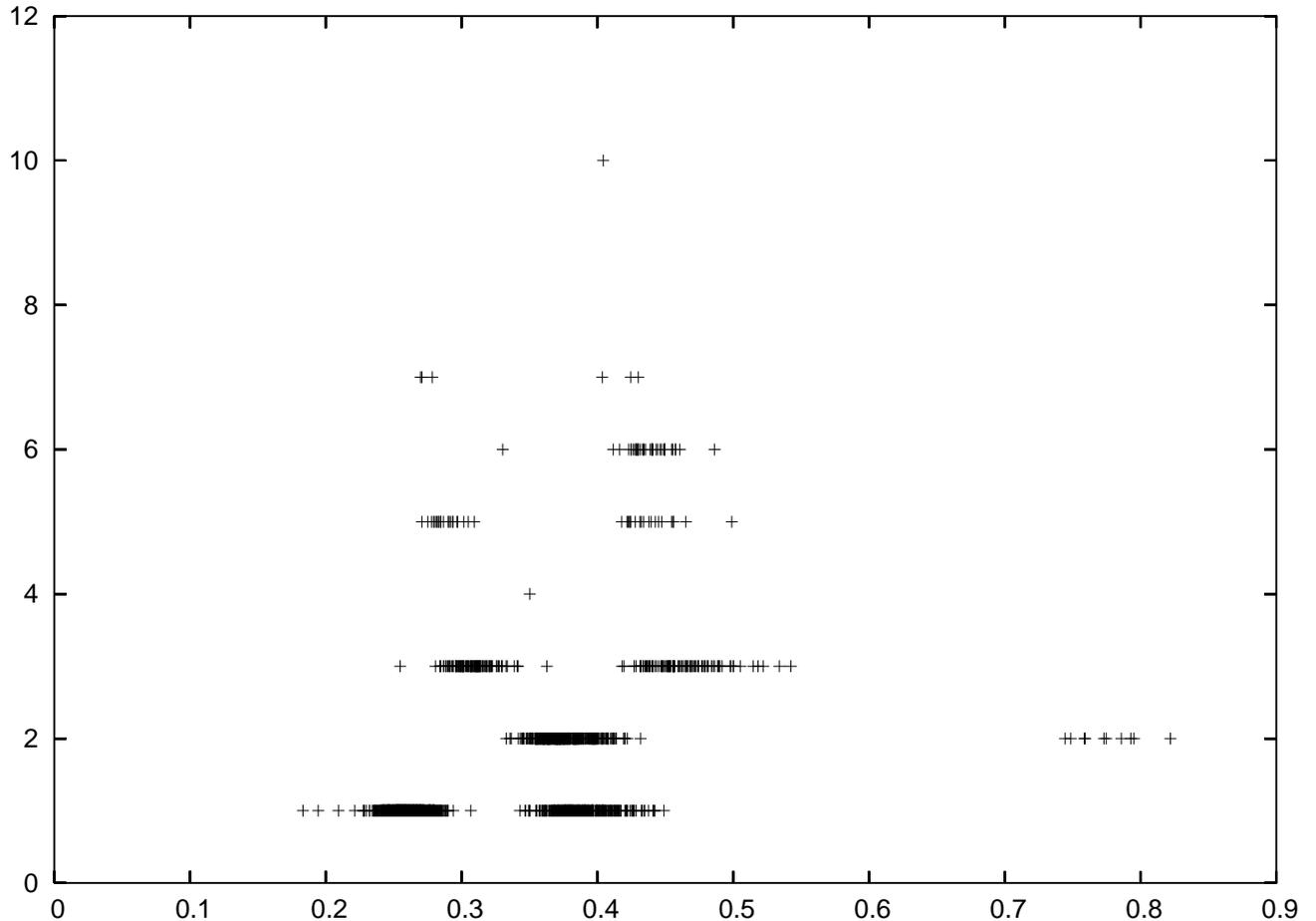,width=5in,angle=-90}
    }
    \caption
    {The rightmost values of the graphs in \ref{fig:c_E all curves} plotted
     against the order of the torsion subgroups of the corresponding elliptic
     curves. We see, for example, that the smallest constants tend to go along
     with the curves whose torsion is trivial. The few curves at the far right with torsion
     size equal to 2 have $c(|d|)$ divisible by 2 when $|d|\neq 4$, and hence,
     because the corresponding discretization is twice as large, the $c(|d|)$'s are four times as 
     likely to vanish.}
    \label{fig:torsion}
\end{figure}

Next we sort the graphs in Figure \ref{fig:c_E all curves} by their
rightmost values at $X=10^8$, and, in Figure \ref{fig:torsion},
plot these values against
the order of the torsion subgroup of the corresponding elliptic curves.
The plot shows that curves with trivial torsion tend to have smaller constants,
followed by curves of torsion size equal to 7 or 5 or 3, then 2 and 1 again, etc.
We do not yet have an explanation for this phenomenon and do not
know how to incorporate it into our model. It seems\cite{kn:rub04} that
to nail down the constant $c^{\pm}_E$ one would need to incorporate
into the model Delaunay's heuristics for Tate-Shavarevich groups \cite{kn:delaunay}. However,
it appears that for primes $p$ dividing the order of the torsion subgroup, the probability
that $c(|d|)$ is divisible by $p$ deviates in a way we do not yet 
understand from a prediction made by
Delaunay for the probability  that the order of the Tate-Shavarevich group is
divisible by a given prime. Delaunay's predictions are for all elliptic curves
sorted by conductor and here we are examining a skinny set of elliptic curves,
namely quadratic twists of a fixed elliptic curve.

\begin{table}[h]
\centerline{\small
\begin{tabular}{|l|l|l|l|l|}
\hline
$11A_i$ &  [0, -1, 1, -10, -20]&  $2.91763323388$&  $13931691$ & 2 \\ \hline  
\multicolumn{5}{|l|}{1, 3, 4, 5, 9 mod 11} \\ \hline
\multicolumn{5}{|l|}{[1/2, -1/2]}\\\hline
\multicolumn{5}{|p{7in}|}{[3, 15, 15, -14, -2, -2],  [4, 11, 12, 0, -4, 0]} \\ \hline\hline
$11A_r$ & $[0, -1, 1, -10, -20]$& $6.34604652140$ & $13931649$ & 6 \\ \hline
\multicolumn{5}{|l|}{1, 3, 4, 5, 9 mod 11} \\ \hline
\multicolumn{5}{|l|}{[1/10, -1/10, 3/10, -3/10, -2/10, 2/10]} \\ \hline 
\multicolumn{5}{|p{7in}|}{[1, 44, 132, -44, 0, 0],  [4, 12, 121, 0, 0, -4],  [4, 33, 45, -22, -4, 0],  [5, 9, 124, -8, -4, -2],  [5, 36, 36, 28, 4, 4],  [16, 16, 25, -4, -4, -12]} \\ \hline\hline
$14A_i$ &  [1, 0, 1, 4, -6] & 5.30196495873& 4432803 & 2 \\ \hline
\multicolumn{5}{|l|}{15, 23, 39 mod 56} \\ \hline
\multicolumn{5}{|l|}{[1/4, -1/4]} \\ \hline
\multicolumn{5}{|p{7in}|}{[4, 15, 15, 2, 4, 4],  [7, 8, 16, -8, 0, 0]} \\ \hline\hline
$15A_i$ &[1, 1, 1, -10, -10] & 3.19248444426& 4749434 & 2 \\ \hline
\multicolumn{5}{|l|}{2, 8 mod 15} \\ \hline 
\multicolumn{5}{|l|}{[1/4, -1/4]} \\ \hline
\multicolumn{5}{|p{7in}|}{[3, 20, 20, -20, 0, 0],  [8, 8, 15, 0, 0, -4]} \\ \hline \hline
$17A_i$ &[1, -1, 1, -1, -14]& 2.74573911809& 14353828 & 2 \\ \hline
\multicolumn{5}{|l|}{3, 5, 6, 7, 10, 11, 12, 14 mod 17} \\ \hline
\multicolumn{5}{|l|}{[1/2, -1/2]} \\ \hline 
\multicolumn{5}{|p{7in}|}{[3, 23, 23, -22, -2, -2],  [7, 11, 20, -8, -4, -6]} \\ \hline \hline
$19A_i$ & [0, 1, 1, -9, -15] & 4.12709239172&  14438275&  2 \\ \hline
\multicolumn{5}{|p{7in}|}{1, 4, 5, 6, 7, 9, 11, 16, 17 mod 19} \\ \hline
\multicolumn{5}{|p{7in}|}{[1/2, -1/2]} \\ \hline
\multicolumn{5}{|p{7in}|}{[4, 19, 20, 0, -4, 0],  [7, 11, 23, -10, -6, -2]} \\ \hline \hline
$19A_r$ &  [0, 1, 1, -9, -15]& 4.07927920046& 14438248& 12 \\ \hline
\multicolumn{5}{|p{7in}|}{1, 4, 5, 6, 7, 9, 11, 16, 17 mod 19}  \\ \hline
\multicolumn{5}{|p{7in}|}{[-1/6, 1/6, -1/6, 1/6, -1/3, -1/3, -1/3, 1/3, 1/3, 1/3, 1/3, -1/3]}  \\ \hline
\multicolumn{5}{|p{7in}|}{[1, 76, 380, -76, 0, 0],  [4, 20, 361, 0, 0, -4],  [4, 77, 96, 40, 4, 4],  [5, 16, 365, 16, 2, 4],  [5, 61, 92, 16, 4, 2],  [5, 76, 92, -76, -4, 0],  [9, 44, 77, 28, 6, 8],  [16, 24, 77, 20, 8, 4],  [17, 44, 44, 12, 16, 16],  [20, 24, 73, 4, 8, 20],  [20, 36, 45, 20, 16, 12],  [25, 36, 36, -4, -16, -16]} \\ \hline
\hline
$21A_i$ & [1, 0, 0, -4, -1]  &3.82197956150 & 4986931 & 2 \\ \hline
\multicolumn{5}{|p{7in}|}{10, 13, 19 mod 21}  \\ \hline
\multicolumn{5}{|p{7in}|}{[1/4, -1/4]}  \\ \hline
\multicolumn{5}{|p{7in}|}{[3, 28, 28, -28, 0, 0],  [7, 12, 24, -12, 0, 0]} \\ \hline
\hline
$26A_i$ & [1, 0, 1, -5, -8] &  3.47934348343 & 4704178 & 2 \\ \hline
\multicolumn{5}{|p{7in}|}{7, 15, 31, 47, 63, 71 mod 104}  \\ \hline
\multicolumn{5}{|p{7in}|}{[-1/2, 1/2]}  \\ \hline
\multicolumn{5}{|p{7in}|}{[7, 15, 31, -14, -6, -2],  [8, 15, 28, 4, 8, 8]} \\ \hline
\hline
$26B_i$ & [1, -1, 1, -3, 3]  & 1.80405719338 & 4704185 & 2 \\ \hline
\multicolumn{5}{|p{7in}|}{3, 27, 35, 43, 51, 75 mod 104}  \\ \hline
\multicolumn{5}{|p{7in}|}{[-1/2, 1/2]}  \\ \hline
\multicolumn{5}{|p{7in}|}{[3, 35, 35, -34, -2, -2],  [4, 27, 27, 2, 4, 4]} \\ \hline
\hline
$30A_i$ & [1, 0, 1, 1, 2]  & 18.5342737810 & 1583137 & 2 \\ \hline
\multicolumn{5}{|p{7in}|}{31, 79 mod 120}  \\ \hline
\multicolumn{5}{|p{7in}|}{[1/8, -1/8]}  \\ \hline
\multicolumn{5}{|p{7in}|}{[4, 31, 31, 2, 4, 4],  [15, 16, 16, -8, 0, 0]} \\ \hline
\end{tabular}
}
\caption{Some of the ternary forms from \cite{kn:rt04}.}
\label{tab:ternary1}
\end{table}

\begin{table}[h]
\centerline{\small   
\begin{tabular}{|l|l|l|l|l|}
\hline
$33A_i$ & [1, 1, 0, -11, 0]  & 2.74463335747 & 5224398 & 3 \\ \hline
\multicolumn{5}{|p{7in}|}{5, 14, 20, 23, 26 mod 33}  \\ \hline
\multicolumn{5}{|p{7in}|}{[1/4, 1/4, -1/2]}  \\ \hline
\multicolumn{5}{|p{7in}|}{[3, 44, 44, -44, 0, 0],  [11, 12, 36, -12, 0, 0],  [15, 20, 20, -4, -12, -12]} \\ \hline
\hline
$34A_i$ & [1, 0, 0, -3, 1] &  7.45670022989 & 4784604 & 2 \\ \hline
\multicolumn{5}{|p{7in}|}{19, 35, 43, 59, 67, 83, 115, 123 mod 136}  \\ \hline
\multicolumn{5}{|p{7in}|}{[1/4, -1/4]}  \\ \hline
\multicolumn{5}{|p{7in}|}{[4, 35, 35, 2, 4, 4],  [8, 19, 36, 4, 8, 8]} \\ \hline
\hline
$35A_i$ & [0, 1, 1, 9, 1] &  2.20504427610 & 5540991 & 2 \\ \hline
\multicolumn{5}{|p{7in}|}{1, 4, 9, 11, 16, 29 mod 35}  \\ \hline
\multicolumn{5}{|p{7in}|}{[-1/2, 1/2]} \\ \hline
\multicolumn{5}{|p{7in}|}{[4, 35, 36, 0, -4, 0],  [11, 15, 39, -10, -6, -10]} \\ \hline
\hline
$37B_i$ & [0, 1, 1, -23, -50] & 7.07044268094 & 14798224 & 2 \\ \hline
\multicolumn{5}{|p{7in}|}{2, 5, 6, 8, 13, 14, 15, 17, 18, 19, 20, 22, 23, 24, 29, 31, 32, 35 mod 37}  \\ \hline
\multicolumn{5}{|p{7in}|}{[1/2, -1/2]}  \\ \hline
\multicolumn{5}{|p{7in}|}{[8, 19, 39, 2, 8, 4],  [15, 20, 23, -8, -14, -4]} \\ \hline
\hline
$43A_r$ & [0, 1, 1, 0, 0] & 10.9373790599 & 14852746 & 22 \\ \hline
\multicolumn{5}{|p{7in}|}{2, 3, 5, 7, 8, 12, 18, 19, 20, 22, 26, 27, 28, 29, 30, 32, 33, 34, 37, 39, 42 mod 43}  \\ \hline
\multicolumn{5}{|p{7in}|}{[1/4, 1/4, -1/4, -1/4, 1/4, 1/4, -1/4, 1/4, -1/4, -1/4, -1/4, -1/4, -1/4, 1/4, 1/4, 1/4, -1/4, 1/4, 1/4, 1/4, -1/4, -1/4]} \\ \hline
\multicolumn{5}{|p{7in}|}{[5, 69, 929, -34, -2, -2],  [5, 241, 276, 104, 4, 2],  [8, 65, 624, 44, 4, 4],  [8, 108, 389, 88, 8, 4],  [12, 29, 932, -28, -4, -4],  [12, 73, 373, 30, 4, 8],  [20, 29, 621, 22, 12, 16],  [20, 61, 277, 54, 12, 8],  [28, 32, 377, 28, 16, 12],  [28, 33, 376, 28, 12, 16],  [29, 77, 157, 14, 18, 26],  [29, 77, 161, 46, 10, 26],  [32, 69, 157, 26, 12, 24],  [37, 45, 237, -2, -26, -34],  [37, 89, 104, 44, 16, 10],  [37, 93, 104, -60, -16, -2],  [45, 77, 113, 62, 42, 10],  [48, 77, 104, -56, -4, -32],  [48, 77, 108, -48, -28, -32],  [61, 80, 89, 76, 42, 16],  [69, 76, 77, -4, -58, -32],  [69, 77, 80, -44, -8, -58]} \\ \hline
\hline
$67A_i$ & [0, 1, 1, -12, -21] & 6.05993680291 & 14974655 & 3 \\ \hline
\multicolumn{5}{|p{7in}|}{1, 4, 6, 9, 10, 14, 15, 16, 17, 19, 21, 22, 23, 24, 25, 26, 29, 33, 35, 36, 37, 39, 40, 47, 49, 54, 55, 56, 59, 60, 62, 64, 65 mod 67}  \\ \hline
\multicolumn{5}{|p{7in}|}{[1/2, -1, 1/2]}  \\ \hline
\multicolumn{5}{|p{7in}|}{[4, 67, 68, 0, -4, 0],  [15, 36, 39, -16, -14, -4],  [16, 19, 71, 6, 16, 12]} \\ \hline
\hline
$67A_r$ & [0, 1, 1, -12, -21] & 1.27377003655 & 14974660 & 70 \\ \hline
\multicolumn{5}{|p{7in}|}{1, 4, 6, 9, 10, 14, 15, 16, 17, 19, 21, 22, 23, 24, 25, 26, 29, 33, 35, 36, 37, 39, 40, 47, 49, 54, 55, 56, 59, 60, 62, 64, 65 mod 67}  \\ \hline
\multicolumn{5}{|p{7in}|}{[-1/2, 1/2, 1/2, -1, -1/2, -1, -1, 1, -1, -1, -1, -1, 1, 1, -1, 1, -1, -1, 1, -1, -1, 1, 1, -1, 1, 1, -1, -1, 1, 1, -1, -1, 1, -1, -1, 1, 1, -1, -1, 1, 1, -1, 1, 1, 1, -1, 1, 1, 1, 1, -1, 1, -1, 1, -1, 1, 1, -1, 1, 1, -1, 1, -1, -1, -1, -1, 1, 1, -1, 1]} \\ \hline
\multicolumn{5}{|p{7in}|}{[1, 268, 4556, -268, 0, 0],  [4, 68, 4489, 0, 0, -4],  [4, 269, 1140, 136, 4, 4],  [9, 269, 508, 92, 8, 6],  [16, 17, 4493, 2, 16, 4],  [17, 237, 332, -172, -8, -6],  [17, 268, 332, -268, -8, 0],  [21, 77, 753, -38, -2, -6],  [21, 217, 308, 192, 12, 2],  [21, 220, 308, 196, 12, 16],  [24, 157, 328, 48, 20, 8],  [24, 236, 261, -196, -20, -12],  [25, 248, 248, 228, 12, 12],  [29, 56, 760, -44, -16, -8],  [29, 89, 504, 20, 24, 26],  [29, 205, 216, 88, 20, 14],  [29, 205, 224, 124, 16, 14],  [29, 224, 224, 180, 16, 16],  [33, 144, 292, -92, -32, -28],  [33, 173, 236, 124, 8, 18],  [33, 173, 237, -118, -14, -18],  [36, 68, 509, 48, 20, 12],  [36, 173, 216, 120, 4, 16],  [37, 77, 464, 68, 8, 26],  [37, 116, 293, -64, -22, -4],  [37, 188, 189, 84, 10, 32],  [37, 189, 189, 110, 10, 10],  [40, 77, 449, 2, 4, 40],  [40, 77, 457, 18, 36, 40],  [40, 149, 237, -130, -20, -16],  [49, 104, 260, -92, -24, -4],  [60, 77, 277, -62, -4, -16],  [60, 84, 277, 64, 4, 44],  [60, 89, 285, -42, -44, -56],  [64, 129, 160, 4, 36, 44],  [65, 68, 301, -20, -2, -40],  [65, 68, 304, -12, -28, -40],  [65, 140, 140, 12, 32, 32],  [65, 140, 160, -84, -44, -32],  [65, 149, 149, 30, 58, 58],  [65, 149, 157, 110, 34, 58],  [68, 77, 237, 10, 12, 24],  [68, 84, 265, 8, 16, 68],  [68, 96, 217, -28, -8, -60],  [68, 96, 237, 92, 12, 60],  [68, 140, 149, 72, 44, 52],  [73, 77, 240, 32, 20, 46],  [73, 92, 181, 12, 18, 8],  [73, 132, 132, -4, -32, -32],  [76, 77, 272, 48, 12, 72],  [76, 93, 173, -6, -8, -20],  [76, 93, 188, 16, 68, 20],  [77, 84, 193, 40, 22, 12],  [77, 89, 220, -8, -36, -70],  [77, 129, 157, 70, 74, 66],  [77, 132, 148, 92, 52, 56],  [77, 132, 157, 100, 74, 56],  [81, 93, 196, 4, 16, 74],  [84, 116, 153, 88, 32, 60],  [84, 116, 169, -48, -80, -60],  [88, 92, 173, 76, 44, 36],  [88, 96, 173, 80, 44, 52],  [89, 132, 132, -4, -72, -72],  [89, 132, 148, -92, -48, -72],  [92, 92, 181, -12, -12, -84],  [93, 93, 173, -6, -6, -82],  [93, 93, 193, -46, -46, -82],  [96, 96, 157, -16, -16, -76],  [100, 104, 149, 64, 96, 44],  [104, 104, 149, -64, -64, -60]} \\ \hline
\end{tabular}
}
\label{tab:ternary2}
\end{table}

\begin{table}[h]
\centerline{\small   
\begin{tabular}{|l|l|l|l|l|}
\hline
$73A_i$ & [1, -1, 0, 4, -3] & 2.79278430294 & 14992818 & 4 \\ \hline
\multicolumn{5}{|p{7in}|}{5, 7, 10, 11, 13, 14, 15, 17, 20, 21, 22, 26, 28, 29, 30, 31, 33, 34, 39, 40, 42, 43, 44, 45, 47, 51, 52, 53, 56, 58, 59, 60, 62, 63, 66, 68 mod 73}  \\ \hline
\multicolumn{5}{|p{7in}|}{[-1/2, 1/2, -1, 1]} \\ \hline
\multicolumn{5}{|p{7in}|}{[7, 43, 84, -40, -4, -6],  [11, 28, 80, 28, 4, 8],  [15, 39, 40, 20, 8, 2],  [20, 31, 44, -28, -4, -12]} \\ \hline
\hline
$79A_r$ & [1, 1, 1, -2, 0] &11.9016007052 &15008174 &40 \\ \hline
\multicolumn{5}{|p{7in}|}{3, 6, 7, 12, 14, 15, 17, 24, 27, 28, 29, 30, 33, 34, 35, 37, 39, 41, 43, 47, 48, 53, 54, 56, 57, 58, 59, 60, 61, 63, 66, 68, 69, 70, 71, 74, 75, 77, 78 mod 79}  \\ \hline
\multicolumn{5}{|p{7in}|}{[-1/2, -1/2, 1/2, 1/2, -1/2, -1/2, -1/2, 1/2, -1/2, -1/2, -1/2, 1/2, 1/2, -1/2, 1/2, 1/2, 1/2, -1/2, 1/2, -1/2, -1/2, -1/2, -1/2, -1/2, 1/2, 1/2, 1/2, 1/2, -1/2, 1/2, 1/2, -1/2, 1/2, 1/2, 1/2, -1/2, -1/2, 1/2, -1/2, 1/2]} \\ \hline
\multicolumn{5}{|p{7in}|}{[12, 369, 449, 54, 8, 4],  [17, 133, 896, -12, -16, -14],  [17, 224, 521, 20, 6, 8],  [17, 264, 449, -12, -14, -16],  [24, 185, 461, 106, 4, 8],  [28, 192, 385, -112, -12, -4],  [29, 33, 2088, -8, -28, -6],  [33, 116, 544, 92, 32, 12],  [37, 153, 376, 8, 36, 34],  [41, 56, 1044, 4, 24, 40],  [41, 116, 425, -4, -38, -8],  [41, 185, 265, -38, -22, -2],  [41, 185, 276, -84, -32, -2],  [41, 216, 236, 20, 28, 36],  [41, 217, 229, -18, -30, -14],  [48, 112, 377, 28, 40, 4],  [48, 193, 224, -32, -36, -20],  [53, 201, 233, -182, -10, -30],  [56, 96, 377, -44, -32, -4],  [60, 157, 216, -20, -4, -32],  [61, 96, 372, -4, -32, -44],  [68, 113, 284, 104, 12, 20],  [68, 172, 185, -92, -32, -4],  [69, 137, 217, 10, 22, 34],  [69, 161, 193, -74, -38, -14],  [77, 113, 249, -62, -46, -26],  [77, 137, 201, -38, -50, -22],  [85, 113, 233, 82, 46, 42],  [93, 145, 161, 70, 22, 46],  [96, 136, 197, -132, -60, -20],  [96, 140, 185, -128, -16, -44],  [108, 120, 161, -52, -36, -4],  [108, 145, 161, 70, 36, 104],  [112, 116, 185, 100, 76, 12],  [113, 116, 185, -100, -14, -56],  [113, 161, 165, 138, 2, 110],  [116, 137, 185, 134, 100, 96],  [132, 137, 149, 98, 64, 108],  [132, 153, 156, -44, -124, -100],  [140, 145, 153, 134, 132, 64]} \\ \hline
\hline
$89B_i$ & [1, 1, 0, 4, 5] &2.18489393577 &15029294 &7 \\ \hline
\multicolumn{5}{|p{7in}|}{3, 6, 7, 12, 13, 14, 15, 19, 23, 24, 26, 27, 28, 29, 30, 31, 33, 35, 37, 38, 41, 43, 46, 48, 51, 52, 54, 56, 58, 59, 60, 61, 62, 63, 65, 66, 70, 74, 75, 76, 77, 82, 83, 86 mod 89}  \\ \hline
\multicolumn{5}{|p{7in}|}{[1/2, -1/2, 1/2, 1/2, 1/2, -1/2, -1]}  \\ \hline
\multicolumn{5}{|p{7in}|}{[3, 119, 119, -118, -2, -2],  [7, 51, 103, -50, -6, -2],  [12, 31, 92, 4, 12, 8],  [15, 24, 95, 24, 2, 4],  [15, 27, 96, -20, -8, -14],  [19, 23, 95, -14, -10, -18],  [23, 31, 48, 16, 12, 2]} \\ \hline
\hline
$109A_i$ & [1, -1, 0, -8, -7] &5.94280424076 &15060017 &3 \\ \hline
\multicolumn{5}{|p{7in}|}{2, 6, 8, 10, 11, 13, 14, 17, 18, 19, 23, 24, 30, 32, 33, 37, 39, 40, 41, 42, 44, 47, 50, 51, 52, 53, 54, 55, 56, 57, 58, 59, 62, 65, 67, 68, 69, 70, 72, 76, 77, 79, 85, 86, 90, 91, 92, 95, 96, 98, 99, 101, 103, 107 mod 109}  \\ \hline
\multicolumn{5}{|p{7in}|}{[1/2, 1/2, -1]}  \\ \hline
\multicolumn{5}{|p{7in}|}{[11, 40, 119, 40, 2, 4],  [19, 23, 119, -22, -18, -2],  [24, 39, 56, 4, 12, 16]} \\ \hline
\hline
$113A_i$ & [1, 1, 1, 3, -4] &2.85781203904 15064917 &7 \\ \hline
\multicolumn{5}{|p{7in}|}{3, 5, 6, 10, 12, 17, 19, 20, 21, 23, 24, 27, 29, 33, 34, 35, 37, 38, 39, 40, 42, 43, 45, 46, 47, 48, 54, 55, 58, 59, 65, 66, 67, 68, 70, 71, 73, 74, 75, 76, 78, 79, 80, 84, 86, 89, 90, 92, 93, 94, 96, 101, 103, 107, 108, 110 mod 113}  \\ \hline
\multicolumn{5}{|p{7in}|}{[-1/2, 1/2, -1/2, -1, -1/2, 1, 1]} \\ \hline
\multicolumn{5}{|p{7in}|}{[3, 151, 151, -150, -2, -2],  [12, 39, 116, 4, 12, 8],  [19, 24, 119, 24, 2, 4],  [20, 47, 68, -44, -4, -12],  [23, 24, 119, 24, 10, 20],  [23, 40, 59, 20, 2, 8],  [35, 39, 47, -10, -34, -6]}\\ \hline
\hline
$139A_i$ & [1, 1, 0, -3, -4] &5.80133204474 &15089693 &5\\ \hline
\multicolumn{5}{|p{7in}|}{1, 4, 5, 6, 7, 9, 11, 13, 16, 20, 24, 25, 28, 29, 30, 31, 34, 35, 36, 37, 38, 41, 42, 44, 45, 46, 47, 49, 51, 52, 54, 55, 57, 63, 64, 65, 66, 67, 69, 71, 77, 78, 79, 80, 81, 83, 86, 89, 91, 96, 99, 100, 106, 107, 112, 113, 116, 117, 118, 120, 121, 122, 124, 125, 127, 129, 131, 136, 137 mod 139} \\ \hline
\multicolumn{5}{|p{7in}|}{[-1/2, 1/2, 1/2, 1/2, -1]} \\ \hline
\multicolumn{5}{|p{7in}|}{[7, 80, 159, 80, 2, 4],  [11, 52, 152, 52, 4, 8],  [20, 28, 139, 0, 0, -4],  [20, 31, 144, 8, 20, 16],  [24, 47, 71, 2, 12, 8]}\\ \hline
\end{tabular}
}
\label{tab:ternary3}
\end{table}

\begin{table}[h]
\centerline{\small   
\begin{tabular}{|l|l|l|l|l|}
\hline
$179A_i$ & [0, 0, 1, -1, -1] &5.10909732904 &15113724 &7\\ \hline
\multicolumn{5}{|p{7in}|}{1, 3, 4, 5, 9, 12, 13, 14, 15, 16, 17, 19, 20, 22, 25, 27, 29, 31, 36, 39, 42, 43, 45, 46, 47, 48, 49, 51, 52, 56, 57, 59, 60, 61, 64, 65, 66, 67, 68, 70, 74, 75, 76, 77, 80, 81, 82, 83, 85, 87, 88, 89, 93, 95, 100, 101, 106, 107, 108, 110, 116, 117, 121, 124, 125, 126, 129, 135, 138, 139, 141, 142, 144, 145, 146, 147, 149, 151, 153, 155, 156, 158, 161, 168, 169, 171, 172, 173, 177 mod 179} \\ \hline
\multicolumn{5}{|p{7in}|}{[1/2, 1/2, 1/2, 1/2, -1/2, -1/2, -1]} \\ \hline
\multicolumn{5}{|p{7in}|}{[4, 179, 180, 0, -4, 0],  [15, 48, 191, 48, 2, 4],  [15, 51, 192, -44, -8, -14],  [16, 47, 183, 6, 16, 12],  [19, 39, 191, -34, -14, -10],  [20, 39, 184, 8, 20, 16],  [39, 56, 76, -52, -20, -12]}\\ \hline
\hline
$233A_i$ & [1, 0, 1, 0, 11] &1.63933561519 &15133226 &13\\ \hline
\multicolumn{5}{|p{7in}|}{3, 5, 6, 10, 11, 12, 17, 20, 21, 22, 24, 27, 34, 35, 39, 40, 41, 42, 43, 44, 45, 47, 48, 53, 54, 57, 59, 61, 65, 67, 68, 69, 70, 73, 75, 77, 78, 79, 80, 82, 83, 84, 86, 87, 88, 90, 93, 94, 95, 96, 97, 99, 103, 106, 108, 111, 114, 115, 118, 119, 122, 125, 127, 130, 134, 136, 137, 138, 139, 140, 143, 145, 146, 147, 149, 150, 151, 153, 154, 155, 156, 158, 160, 163, 164, 165, 166, 168, 172, 174, 176, 179, 180, 185, 186, 188, 189, 190, 191, 192, 193, 194, 198, 199, 206, 209, 211, 212, 213, 216, 221, 222, 223, 227, 228, 230 mod 233} \\ \hline
\multicolumn{5}{|p{7in}|}{[1/2, -1/2, 1/2, 1, 1/2, 3/2, -3/2, 1, -1, -3, 1, 1, -1]} \\ \hline
\multicolumn{5}{|p{7in}|}{[3, 311, 311, -310, -2, -2],  [11, 87, 255, -82, -6, -10],  [12, 79, 236, 4, 12, 8],  [20, 95, 140, -92, -4, -12],  [24, 39, 239, 2, 24, 4],  [24, 43, 239, 10, 24, 20],  [27, 39, 244, -28, -16, -22],  [35, 80, 84, 28, 24, 4],  [39, 80, 96, -68, -8, -36],  [40, 47, 119, 2, 20, 8],  [44, 68, 87, -36, -20, -28],  [47, 68, 87, 36, 38, 40],  [48, 59, 79, 2, 16, 12]}\\ \hline
\end{tabular}
}
\label{tab:ternary4}
\end{table}

\subsubsection{Acknowledgements}

We wish to thank Gonzalo Tornaria and Fernando Rodriguez-Villegas
for supplying us with a database of ternary quadratic forms which were used to compute
the $L$-values described in this paper. Andrew Granville and William Stein allowed us to
use their computer clusters to perform some of the computations described
here. Atul Pokharel assisted with the preperation of
Figures~\ref{fig:1st vs 2nd}--~\ref{fig:torsion}.
The authors are grateful to AIM and the Isaac Newton Institute for
very generous support and hospitality, as well as to support in the
form of an NSF FRG grant that helped make some of this work possible.
JPK thanks the EPSRC for support in the form of a Senior Research Fellowship,
NCS the Royal Society and EPSRC for support during the
course of this work, MOR the NSF and NSERC, and JBC the NSF.

\end{document}